\documentclass{siamltex}

\usepackage{microtype}
\usepackage[text={5.125in,8.25in},centering]{geometry}
\usepackage{amsmath}
\usepackage{amssymb}
\usepackage{verbatim} 
\usepackage{graphicx}
\usepackage{url}
\usepackage{booktabs}
\usepackage[colorlinks,citecolor={blue},urlcolor={blue},linkcolor={blue}]{hyperref}
\usepackage{gitinfo}

\input{./macros}
\newcommand{\eval}{\mathbb{E}}
\newcommand{\batch}{\Bscr}
\newcommand{\batchcomp}{\Nscr}
\newcommand{\dist}{\pi}

\usepackage{boxedminipage}
\usepackage{color}
\definecolor{softblue}{rgb}{0.90,0.92,1.00}
\makeatletter\newenvironment{btheorem}{%
\begin{lrbox}{\@tempboxa}\begin{minipage}{0.97\textwidth}\begin{theorem}}%
{\end{theorem}\end{minipage}\end{lrbox}%
\par\hbox{}\noindent%
{\setlength{\fboxsep}{0pt}\colorbox{softblue}{\setlength{\fboxsep}{4pt}\begin{boxedminipage}{\textwidth}\usebox{\@tempboxa}\end{boxedminipage}}}%
\vspace{0.5\baselineskip}}
\makeatother
\makeatletter\newenvironment{blemma}{%
\begin{lrbox}{\@tempboxa}\begin{minipage}{0.97\textwidth}\begin{lemma}}%
{\end{lemma}\end{minipage}\end{lrbox}%
\par\hbox{}\noindent%
{\setlength{\fboxsep}{0pt}\colorbox{softblue}{\setlength{\fboxsep}{4pt}\begin{boxedminipage}{\textwidth}\usebox{\@tempboxa}\end{boxedminipage}}}%
\vspace{0.5\baselineskip}}
\makeatother
\makeatletter\newenvironment{bcorollary}{%
\begin{lrbox}{\@tempboxa}\begin{minipage}{0.97\textwidth}\begin{corollary}}%
{\end{corollary}\end{minipage}\end{lrbox}%
\par\hbox{}\noindent%
{\setlength{\fboxsep}{0pt}\colorbox{softblue}{\setlength{\fboxsep}{4pt}\begin{boxedminipage}{\textwidth}\usebox{\@tempboxa}\end{boxedminipage}}}%
\vspace{0.5\baselineskip}}
\makeatother

\usepackage{framed}\definecolor{shadecolor}{rgb}{.9,.92,1}

\title{Hybrid Deterministic-Stochastic\\Methods for Data Fitting}
\author{Michael P.~Friedlander\thanks{Department of Computer Science,
    University of British Columbia, Vancouver V6T 1Z4, B.C., Canada
    ({\tt mpf@cs.ubc.ca}). The work of this author was supported by
    NSERC Discovery Grant 312104.%
  } \and Mark Schmidt\thanks{INRIA - SIERRA team, Laboratoire
    d'Informatique de l'Ecole Normale Superieure (INRIA/ENS/CNRS UMR
    8548), 23 avenue d'Italie, 75214 Paris CEDEX 13, France ({\tt
      mark.schmidt@inria.fr}).  The work of this author was supported
    by the CRD DNOISE II NSERC grant, and the SIERRA grant from the
    European Research Council (SIERRA-ERC-239993).}  \hfill\hbox{Date:
    \gitCommitterDate}}
\begin{document}

\maketitle

\begin{abstract}
  Many structured data-fitting applications 
  require the solution of an optimization problem involving a sum over
  a potentially large number of measurements.  Incremental gradient
  algorithms offer inexpensive iterations by sampling a subset of the
  terms in the sum; these methods can make great progress initially,
  but often slow as they approach a solution. In contrast,
  full-gradient methods achieve steady convergence at the expense of
  evaluating the full objective and gradient on each iteration. We
  explore hybrid methods that exhibit the benefits of both
  approaches. Rate-of-convergence analysis shows that by controlling
  the sample size in an incremental-gradient algorithm, it is possible
  to maintain the steady convergence rates of full-gradient methods.
  We detail a practical quasi-Newton implementation based on this
  approach. Numerical experiments illustrate its potential benefits.
\end{abstract}

\thispagestyle{plain}

\section{Introduction}

Data-fitting applications are often typified by optimization problems
of the form
\begin{equation}
  \label{eq:1}
  \minimize{x\in\Real^n}\quad f(x) \defd \frac1M\sum_{i=1}^M f_i(x),
\end{equation}
where each function $f_i$ corresponds to a single observation or
measurement, and models the misfit for a given choice of parameters
$x$. The aim is to choose parameters that minimize the misfit (or loss)
across all measurements. The canonical example is least squares, and
in that case,
$$
f_i(x) = \half(a_i\T x - b_i)^2 \text{for all} i=1,\ldots,M.
$$
This misfit model corresponds to a linear model with Gaussian errors
on the measurements $b$. For applications where the measurements $b$
are binary, a more appropriate model is logistic regression,
described by the choice
$$
f_i(x) = \log(1+\exp[-b_i\drop a_i\T x]) \text{for all} i=1,\ldots,M.
$$
These are both special cases of the more general maximum-likelihood
problem. The maximum-likelihood approach gives rise to separable
problems like~\eqref{eq:1} whenever the measurements $(a_i,b_i)$ are
assumed to be independent and identically distributed.

If the number of measurements $M$ is very large, or if the individual
$f_i$ are complicated functions (e.g., each $f_i$ evaluation may
require the solution of a partial differential equation), then
evaluating $f(x)$ and $\nabla f(x)$ can be computationally expensive.
However, there is often a large amount of uniformity in the
measurements, which means that a full evaluation of $f(x)$ and $\nabla
f(x)$ may be unnecessary to make progress in solving~\eqref{eq:1}.
This motivates \emph{incremental-gradient} methods, in which each
iteration only evaluates the gradient with respect to a single
$f_i$~\cite[\S3.2]{bertsekas1996neuro}.


Incremental-gradient methods enjoy an iteration cost that is $M$ times
faster than full-gradient methods because the iterations are
independent of $M$.  Thus, in the time it takes to make one
full-gradient iteration, the incremental-gradient method can achieve
$M$ iterations, which often results in rapid initial progress.
However, the number of iterations needed to reach the same level of
accuracy may be much
higher. 
Indeed, because of their faster convergence rate, full-gradient methods
must eventually dominate incremental-gradient methods.

Our aim is to develop a method that exhibits the benefits of these two
extremes. The approach is based on starting with iterations that
resemble an incremental-gradient approach and use relatively few
measurements to approximate the gradient; as the iterations proceed,
the algorithm gradually increases the number of measurements. This
preserves the rapid initial progress of incremental-gradient methods
without sacrificing the convergence rate of full-gradient methods.

\subsection{Gradient descent with error}
\label{sec:gradient-descent0-with-error}

In the most basic version of the algorithm that we consider, each
iterate is computed via the update
\begin{equation}\label{eq:2}
x\kp1 = x\k - \alpha\k g\k
\end{equation}
for some step size $\alpha\k$; the search direction
\begin{equation}\label{eq:3}
g\k \defd \nabla f(x\k) + e\k
\end{equation}
is an approximation of the gradient, and $e\k$ is the residual in its
computation. 
This framework is the point
of departure in the analysis of
\S\S\ref{sec:convergence}--\ref{sec:batch-error}.
Evidently, the full gradient method (i.e., steepest
descent) corresponds to the case where $e\k\equiv0$, which means that
the gradient is exactly computed at each iteration. 

The stochastic approximation method deals with the case where the
residual $e\k$ is a random variable; see~\cite{bertsekas1996neuro}. In
this context we typically assume that $\mathbb{E}[e\k] = 0$ and
$\mathbb{E}[\norm{e\k}^2]\le B$ for some constant $B$, which are
sufficient to guarantee that the iterates converge in a probabilistic
sense (for a suitable choice of decreasing step sizes).
Incremental-gradient methods are a special case of stochastic
approximation.  Here, rather than computing the full gradient $\nabla
f(x\k)$ on each iteration, a function $f_i$ is randomly selected among
$i\in\{1,\ldots,M\}$, and the gradient estimate is constructed as $g\k
= \nabla f_i(x\k)$.

In this work we consider scenarios in which the gradient residual can
be controlled on each iteration, by specifying a per-iteration bound
$B_k$ on the norm of the residual. In particular, we characterize
convergence of the gradient-with-error
algorithm~\eqref{eq:2}--\eqref{eq:3} under two different conditions:
\begin{equation}
  \label{eq:8}
    \hbox{deterministic:}\quad \norm{e\k}^2 \le B\k;
    \textt{or}
    \hbox{stochastic:}\quad \mathbb{E}[\norm{e\k}^2] \le B\k.
\end{equation}
Our analysis applies whether the noise is deterministic or
stochastic, and we do not assume that the noise has zero mean.  In the
context of problem~\eqref{eq:1}, the error in the gradient is a result
of using a \emph{sample} of the $f_i$ functions (sometimes referred to
as a \emph{batch}).  In~\S\ref{sec:batch-error} we show how the
\emph{sample size} (or \emph{batch size}) influences the bound $B\k$.

We also consider in \S\ref{sec:practicalImplementation} the case in
which the approximate gradient $g\k$ is scaled to account for
curvature information in $f$. In this case, the iteration update is
\begin{equation}
  \label{eq:25}
  x\kp1 = x\k - \alpha\k d\k,
\end{equation}
where $d\k$ solves the system
\begin{equation}
  \label{eq:27}
  H\k d = -g\k,
\end{equation}
and $H\k$ is a positive-definite approximation (e.g., a
quasi-Newton Hessian) to $\nabla^2 f$.

\subsection{Assumptions and notation}
\label{sec:assumptions}
We make the blanket assumption throughout that a minimizer $\xstar$ of
$f$ always exists, that the functions $f_i:\Real^n\to\Real$ are
continuously differentiable, and that the overall gradient of $f$ is
uniformly Lipschitz continuous, i.e., for some positive $L$,
\begin{subequations} \label{eq:assumptions}
\begin{equation}\label{eq:12}
\norm{\nabla f(x) - \nabla f(y)} \le L \norm{x-y} \quad\hbox{for
  all}\quad x,y\in\Real^n.
\end{equation}
We also assume that $f$ is strongly convex with
(positive) parameter $\mu$:
\begin{equation} \label{eq:4}
f(y) \ge f(x) + (y-x)^T\nabla f(x) + \half \mu \norm{y-x}^2
 \quad\hbox{for all}\quad x,y\in\Real^n.
\end{equation}
\end{subequations}
If $f$ is twice-continuously differentiable, then these
assumptions are equivalent to the condition that the eigenvalues of
the Hessian are uniformly bounded above and below:
\[
 \mu I \preceq \nabla^2 f(x) \preceq L I.
\]
The ratio $L/\mu\ge1$ is known as the \emph{condition} number of $f$
\cite[\S2.1.3]{nesterov2004introductory}. 

We make the assumption---standard in the stochastic optimization
literature~\cite[\S4.2]{bertsekas1996neuro}---that
\begin{equation}
  \label{eq:incGradBnd}
  \norm{\nabla f_i(x)}^2 \leq \beta_1 + \beta_2\norm{\nabla f(x)}^2
  \quad\hbox{for all $x$ and $i=1,\ldots,M$},
\end{equation}
for some constants $\beta_1\geq 0$ and $\beta_2\geq 1$.  This implies that
$\norm{e\k}$ is bounded as a function of the true gradient of the objective.

We describe two versions of linear convergence of function values.
The first is denoted as \emph{weak} linear convergence, and is
characterized by a bounding sequence on the function value at every
iteration $k$:
\begin{equation} \label{eq:weak-rate}
f(x\k) - f(\xstar) = \Oscr(\sigma^k) \text{for some} \sigma<1.
\end{equation}
This is a non-asymptotic version of R-linear
convergence. (See, for example, \cite[\S A.2]{NoceWrig:2006}.) The
second is denoted as \emph{strong} linear convergence, and
characterizes the decrease of the function value at every iteration
$k$:
\begin{equation} \label{eq:strong-rate}
f(x\kp1) - f(\xstar) \leq \sigma [f(x\k) - f(\xstar)] \text{for some} \sigma<1.
\end{equation}
We emphasize that this inequality applies to all iterations of the
algorithm, and is a non-asymptotic version of Q-linear
convergence. Note that \eqref{eq:strong-rate} implies
\eqref{eq:weak-rate}, though the converse is not implied.

All of the convergence results that we analyze are
described only in terms of convergence of the objective function
values, and not the iterates $x\k$. This is sufficient, however,
because the strong convexity assumption allows us to directly deduce a
convergence rate of $x\k$ to $\xstar$ via the corresponding function
values. In particular, strong convexity of $f$ implies
\begin{equation}
\label{eq:xstar_fxstar}
\frac{\mu}{2}\norm{x\k - \xstar}^2 \le f(x\k) - f(\xstar).
\end{equation}
Thus, the rate at which the squared error $\norm{x\k-\xstar}^2$ converges is at
least as fast as the rate at which $f(x\k)$ converges to the optimal
value $f(\xstar)$.

If the matrix $H\k$ in the iterations~\eqref{eq:25}--\eqref{eq:27} is uniformly
positive definite and uniformly bounded in norm (as can be enforced in
practice), then assumptions~\eqref{eq:12}--\eqref{eq:4} can be
replaced by the following conditions: there exist positive constants
$L'$ and $\mu'$ such that for all $x,y\in\Real^n$ and for all
$k=0,1,\ldots$,
\begin{subequations} \label{eq:assumptions-qn}
  \begin{align}
    \norm{\nabla f(x) - \nabla f(y)}_{H\inv\k} \le L' \norm{x-y}_{H\k}\phantom,
    \\ f(y) \ge f(x) + (y-x)^T\nabla f(x) + \half \mu' \norm{x-y}^2_{H\k},
  \end{align}
\end{subequations}
where the quadratic norm $\norm{x}_{H\k}=\sqrt{x\T H\k x}$ and its
dual $\norm{x}_{H\inv\k}$ are used instead of the Euclidean norm. It
can then be verified that all of the results in
\S\S\ref{sec:convergence}--\ref{sec:batch-error} apply to the
Newton-like algorithm~\eqref{eq:25}--\eqref{eq:27}, where the
parameters $L$ and $\mu$ in those results are replaced by the
parameters $L'$ and $\mu'$ found in~\eqref{eq:assumptions-qn}. The
benefit of this approach is that a judicious choice of the scaling
$H\k$ can lead to a scaled condition number $L'/\mu'$ that can be
smaller than the condition number $L/\mu$ of the unscaled objective
$f$, effectively improving on the error constants found in the
convergence results.

\subsection{Contributions} \label{sec:approach}

This paper is divided into six components.
\paragraph{Weak linear convergence with generic bounds \rm(\S\ref{sec:weak-convergence})}
We analyze the convergence rate under a generic sequence $\{B\k\}$.
Our results imply that for any (sub)linearly decreasing sequence
$\{B\k\}$, the algorithm has a weak (sub)linear convergence rate.  In
the expected-error version of~\eqref{eq:8}, the convergence rate is
described in terms of the expected function value.

\paragraph{Strong linear convergence with particular bounds \rm (\S
  \ref{sec:strong-convergence})}
We describe a particular construction of the sequence $\{B\k\}$ that
ensures that the algorithm has a strong linear convergence rate.  The
rate achieved under this sequence can be arbitrarily close to the rate
of the standard gradient method without error, without requiring an
exact gradient calculation on any iteration.

\paragraph{Sublinear convergence without strong convexity \rm
  (\S\ref{sec:weakly})}  
Without the strong convexity assumption on $f$, the convergence rate
for the deterministic gradient method is sublinear. We show that a
summable sequence $\{B\k\}$ is sufficient to maintain the same sublinear rate.

\paragraph{Application to sample-average gradients \rm(\S\ref{sec:batch-error})}
For data-fitting problems of the form~\eqref{eq:1}, we show that a
growing sample-size strategy can be used as a mechanism for
controlling the error in the estimated gradient and achieving a linear
rate.
In effect, choosing the sample size allows us to control the 
error, as in~\eqref{eq:8}.  By growing the sample size sufficiently
fast, we implicitly set the rate at which $B\k\to0$, and hence the
overall rate of the algorithm.


\paragraph{A practical quasi-Newton implementation \rm(\S\ref{sec:practicalImplementation})}
We describe a practical implementation of the ideas based on a
limited-memory quasi-Newton approximation and a heuristic line search. 

\paragraph{Numerical results \rm(\S\ref{sec:experiments})}
We evaluate the implementation on a variety of data-fitting
applications (comparing to incremental gradient methods and a deterministic 
quasi-Newton method).

\subsection{Related work}
\label{sec:related-work}

Our approach is based on bridging the gap between two ends of a
spectrum, where incremental gradient methods are at the end of ``cheap
iterations with slow convergence'', and full gradient methods are
``expensive iterations with fast convergence''.  In particular,
incremental-gradient methods achieve an expected sublinear convergence
rate on the expected value of $f(x\k)$, i.e.,
\[
\mathbb{E}[f(x\k) - f(\xstar)] = \mathcal{O}(1/k),
\]
where the iterations are described by \eqref{eq:2} for
$\alpha\k=\Oscr(1/k)$~\cite[\S2.1]{nemirovski2009robust}. In fact,
among all first-order methods, this is the best possible dependency on
$k$ given only a first-order stochastic oracle; thus a linear rate is not
possible~\cite[\S14.1]{nemirovski1994efficient}.

In contrast, consider the basic gradient-descent iteration \eqref{eq:2} with a
fixed step size $\alpha\k=1/L$ and $g\k = \nabla f(x\k)$. It is
well known that this algorithm has a strong linear convergence rate,
and satisfies the per-iteration decrease in~\eqref{eq:strong-rate}
with $\sigma=1-\mu/L$;
see~\cite[\S8.6]{luenberger2008linear}.

Some authors have analyzed incremental gradient methods with a
constant step size~\cite{nedic2000convergence}.  Although this
strategy does not converge to the optimal solution, it does converge
at a linear rate to a neighborhood of the solution (where the size of
the neighborhood increases with the step size).

In the context of incremental gradient methods, several other hybrid
methods have been proposed that achieve a linear convergence rate. The
works of Bertsekas~\cite{bertsekas1997new} and Blatt, Hero, and
Gauchman~\cite{blatt2008convergent} are the closest in spirit to our
proposed approach.  However, the convergence rates for these methods
treat full passes through the data as iterations, similar to the full
gradient method.  Further, there are numerical difficulties in
evaluating certain sequences associated with the method of Bertsekas,
while the method of Blatt et al.\@ may require an excessive amount of
memory.

This is not the first work to examine a growing sample-size
strategy. This type of strategy appears to be a ``folk'' algorithm
used by practitioners in several application domains, and it is
explicitly mentioned in an informal context by some authors such as
Bertsekas and Tsitsiklis~\cite[page 113]{bertsekas1996neuro}, who
refer to growing sample-size techniques as ``batching'' strategies.
This is the first work, that we are aware of, that presents a
theoretical analysis of the technique and that proposes a practical
large-scale quasi-Newton implementation along with an experimental
evaluation.

Gradient descent with a decreasing sequence of errors in the gradient
measurement was previously analyzed by Luo and
Tseng~\cite{luo1993error}, and they present analogous weak and
strong linear convergence results depending on the sequence of bounds on
the noise.  Our analysis extends this previous work in several ways:
\begin{itemize}
\item The weak linear convergence rate described
  in~\cite{luo1993error} requires a per-iteration strict decrease in
  the objective.
  In contrast, our analysis in~\S\ref{sec:weak-convergence} does not
  require this assumption, and allows for the (realistic) possibility
  that the noisy gradient can lead to an increase in the objective
  function on some iterations.
\item The strong linear convergence rate shown in~\cite{luo1993error}
  only holds asymptotically, while the construction we give
  (see \S\ref{sec:strong-convergence}) leads to a non-asymptotic rate of the form~\eqref{eq:strong-rate}
  that applies to all iterations of the algorithm.
\item Luo and Tseng~\cite{luo1993error} consider deterministic errors
  in the gradient measurement that can be bounded in an absolute
  sense; we also consider the more general scenario where the error is
  stochastic and can only be bounded in expectation.
\end{itemize}

\subsection{Reproducible research} \label{sec:repr-rese}

Following the discipline of reproducible research, the source code and
data files required to reproduce the experimental results of
this paper can be downloaded from
\begin{center}\url{http://www.cs.ubc.ca/labs/scl/FriedlanderSchmidt2011}.\end{center}

\section{Convergence analysis}
\label{sec:convergence}

Our convergence analysis first considers a basic first-order method
with the constant step size $\alpha\k = 1/L$.
The following intermediate result
establishes an upper bound on the objective value at each iteration in
terms of the residual in the computed gradient.
\begin{blemma}\label{le:upper-bound}
At each iteration $k$ of algorithm~\eqref{eq:2}, with $\alpha\k\equiv 1/L$,
\begin{equation}
\label{eq:linCon1}
  f(x\kp1) - f(\xstar) \le
  (1 - \mu/L)[f(x\k) - f(\xstar)] + \frac{1}{2L}\norm{e\k}^2.
\end{equation}
\end{blemma}

\begin{proof}
  It follows from assumptions~\eqref{eq:assumptions} that the
  following inequalities hold:
\begin{subequations}
\begin{align}
  f(y) &\leq f(x) + (y - x)^T\nabla f(x) + \frac{L}{2}\norm{y - x}^2,
  \label{eq:lip}
\\f(y) &\geq f(x) + (y - x)^T\nabla f(x) + \frac{\mu}{2}\norm{y - x}^2.
  \label{eq:sc}
\end{align}
\end{subequations}
Use $x=x\k$ and $y=x\k - (1/L) g\k$ in~\eqref{eq:lip} and
simplify to obtain
\begin{equation*}
  f(x\k - (1/L) g(x\k)) 
  \le f(x\k) - \frac{1}{L}g(x\k)^T\nabla f(x\k) + \frac{1}{2L}\norm{g(x\k)}^2.
\end{equation*}
Next, use the definitions of $x\kp1$ and $g\k$
(cf.~\eqref{eq:2}--\eqref{eq:3}) in this expression to obtain
\begin{equation} \label{eq:6}
\begin{aligned}
  f(x\kp1) & \leq f(x\k) - \frac{1}{L}(\nabla f(x\k)
     + e\k)^T\nabla f(x\k) + \frac{1}{2L}\norm{\nabla f(x\k) + e\k}^2 \\
  & = f(x\k) - \frac{1}{L}\norm{\nabla f(x\k)}^2
  - \frac{1}{L}\nabla f(x\k)^Te\k \\
  &\phantom{=f(x\k)} +\frac{1}{2L}\norm{\nabla f(x\k)}^2
  + \frac{1}{L}\nabla f(x\k)^Te\k + \frac{1}{2L}\norm{e\k}^2\\
  & = f(x\k) - \frac{1}{2L}\norm{\nabla f(x\k)}^2 + \frac{1}{2L}\norm{e\k}^2.
\end{aligned}
\end{equation}

We now use~\eqref{eq:sc} to derive a lower bound on the norm of
$\nabla f(x\k)$ in terms of the optimality of $f(x\k)$. Do this by
minimizing both sides of~\eqref{eq:sc} with respect to $y$: by
definition, the minimum of the left-hand side is achieved by $y=\xstar$;
the minimizer of the right-hand side is given by $y = x -
(1/\mu)\nabla f(x)$.  Thus,
\begin{equation*}
f(\xstar) \geq f(x) - \frac{1}{\mu}\nabla f(x)^T\nabla f(x)
                   + \frac{1}{2\mu}\nabla f(x)^T\nabla f(x)
       = f(x) - \frac{1}{2\mu}\norm{\nabla f(x)}^2
\end{equation*}
for any $x$.  Re-arranging and specializing to the case where $x=x\k$,
\begin{equation} \label{eq:grdBnd}  
  \norm{\nabla f(x\k)}^2 \geq 2\mu[f(x\k) - f(\xstar)].
\end{equation}
Subtract $f(\xstar)$ from both sides of~\eqref{eq:6} and
use~\eqref{eq:grdBnd} to get
\begin{align*} \label{eq:linCon}
  f(x\kp1) - f(\xstar) & \leq f(x\k) - f(\xstar)
  - \frac{\mu}{L}[f(x\k)-f(\xstar)]
  + \frac{1}{2L}\norm{e\k}^2\nonumber\\
  & = (1 - \mu/L)[f(x\k) - f(\xstar)] + \frac{1}{2L}\norm{e\k}^2,
\end{align*}
which gives the required result.
\end{proof}

As an aside, note that~\eqref{eq:6} shows that the objective decreases
monotonically, i.e., $f(x\kp1) < f(x\k)$ if $\norm{e\k} <
\norm{\nabla f(x\k)}$. In general, however, we do not require this
condition.

\subsection{Weak linear convergence}
\label{sec:weak-convergence}

In this section we show that if $\{B\k\}$ is any (sub)linearly
convergent sequence, then algorithm~\eqref{eq:2} has a (sub)linear
convergence rate. This result reflects that the convergence rate of
the approximate gradient algorithm is not better than the rate at
which the noise goes to zero, and of course is also not better than
the rate of the noiseless algorithm.

\begin{btheorem}[Weak convergence rate under absolute error bounds]
  \label{th:linear-rate-weak}
  Suppose that $\norm{e\k}^2\le B\k$, where
  \begin{equation}\label{eq:20}
  \lim_{k\to\infty} B\kp1/B\k \le 1.
  \end{equation}
  Then at each iteration of algorithm~\eqref{eq:2} with
  $\alpha\k\equiv1/L$, for any $\epsilon>0$ we have
  \begin{equation*}
  f(x\k) - f(\xstar) \le (1-\mu/L)^k[f(x_0)-f(x_*)] + \Oscr(C\k),
  \end{equation*}
  where $C\k=\max\{B\k, (1-\mu/L+\epsilon)^k\}$. 
\end{btheorem} 
\begin{proof}
Let $\rho:=1-\mu/L$. Because $\norm{e\k}^2 \leq B\k$, Lemma~\ref{le:upper-bound} implies
\begin{equation*}
  f(x\kp1) - f(\xstar) \le [f(x\k) - f(\xstar)]\rho + \frac{1}{2L}B\k.
\end{equation*}
Applying this recursively,
\begin{equation*}
f(x\k) - f(\xstar) = [f(x_0)-f(x_*)]\rho^k + \Oscr(\mu\k),
\end{equation*}
where
$
\mu\k \defd \sum_{i=0}^{k-1}\rho^{k-i-1}B_i.
$
Observe that $\mu\kp1=B\k + \rho\mu\k$.

We now show that $\mu_k\le\xi C_k$ for all $k$ for some
$\xi>0$. It follows from~\eqref{eq:20} and the definition of $C\k$
that $\lim_{k\to\infty}C_{k+1}/C_k \geq \rho+\epsilon$, and thus 
that there exists some $N$ such that $C_{k+1}/C_{k} - \rho \geq \epsilon/2$ for all
$k \geq N$.  
We now choose $\xi$ such that
\[
  \hbox{$\mu_{k} \leq \xi C_{k}$ for all  $k \leq N$,}
  \text{and}
  \xi \epsilon/2 \geq 1.
\]
This is always possible because the $\mu_{k}$ are finite and $\epsilon > 0$.  
We now show by induction that $\mu_k \leq \xi C_k$ for all $k$.  
This trivially holds for all $k \leq N$ by the definition of $\xi$.
  Assuming this
holds for some arbitrary $k \geq N$, we have
\begin{align*}
\mu_{k+1}
  = B\k + \rho \mu_k
 & \leq C\k + \rho\xi C\k \\
 & \leq \xi(\epsilon/2)C\k + \rho\xi C\k
\leq \xi( C\kp1/C\k  -\rho)C\k + \rho\xi C\k
 = \xi C\kp1.
\end{align*}
Thus, $\mu_k = \Oscr(C\k)$ for all $k$, as required.
\end{proof}

An implication of this result is that the algorithm has a linear convergence
rate if $B\k$ converges linearly to zero.  For example, if $B\k = \Oscr(\gamma^k)$ with
$\gamma < 1$, then
\begin{equation*}
 f(x\k) - f(\xstar) = \Oscr(\sigma^k),
\end{equation*}
where $\sigma=\max\{\gamma,\, (1-\mu/L + \epsilon)\})$ for any
positive $\epsilon<\mu/L$.

Theorem~\ref{th:linear-rate-weak} also yields a convergence rate in scenarios where the
high cost of computing an accurate gradient might make it appealing to
allow the error in the gradient measurement to decrease sublinearly.
For example, if $B\k = \Oscr(1/k^2)$, then $f(x\k)-f(\xstar) =
\Oscr(1/k^2)$, which is the rate achieved by Nesterov's optimal method
for general smooth convex (but not necessarily strongly
convex) functions in the noiseless
setting~\cite[\S2.1]{nesterov2004introductory}.
Theorem~\ref{th:linear-rate-weak} also allows for the possibility
that the bound $B\k$ does not converge to zero (necessarily implying
the limit of $B\kp1/B\k$ is one). In this case, the result simply states that the distance
to optimality is eventually bounded by a constant times $B\k$.

The above analysis allows for the possibility that the approximate gradient is computed
by a stochastic algorithm where the error made by the algorithm can be bounded
in an absolute sense.  We now consider the more general case where the error
can only be bounded in expectation.  The following result is the counterpart to
Theorem~\ref{th:linear-rate-weak}, where we instead have a bound on
the expected value of $\norm{e\k}^2$.

\begin{btheorem}[Weak expected convergence rate under expected error bounds]
  \label{th:linear-rate-weak-expected}
  Suppose that $\eval[\norm{e\k}^2]\le B\k$, where
  \begin{equation*}
  \lim_{k\to\infty} B\kp1/B\k \le 1.
  \end{equation*}
  Then at each iteration of algorithm~\eqref{eq:2} with
  $\alpha\k\equiv1/L$, for any $\epsilon>0$ we have
  \begin{equation*}
  \eval[f(x\k) - f(\xstar)] \le (1-\mu/L)^k[f(x_0)-f(x_*)]+\Oscr(C\k),
  \end{equation*}
  where $C\k=\max\{B\k, (1-\mu/L+\epsilon)^k\}$.
\end{btheorem} 
\begin{proof}
  We use Lemma~\ref{le:upper-bound} and take expectations
  of \eqref{eq:linCon1} to obtain
  \[
  \mathbb{E}[f(x\kp1) - f(\xstar)] \leq (1 - \mu/L)\mathbb{E}[f(x\k) - f(\xstar)]
  + \frac{1}{2L}\mathbb{E}[\norm{e\k}^2].
  \]
  Proceeding as in the proof of Theorem~\ref{th:linear-rate-weak}, we
  obtain a similar result, but based on the expected value of the
  objective.
\end{proof}

\subsection{Strong linear convergence}
\label{sec:strong-convergence}

We now describe a particular construction of the sequence $\{B\k\}$ that
allows us to achieve a linear decrease of the function values that
applies to every iteration $k$.
This strong guarantee, however, comes at the price of requiring bounds
on three quantities that are unknowable for general problems: the
strong convexity constant $\mu$, the gradient's Lipschitz constant
$L$, and a non-trivial lower bound on the current iterate's distance
to optimality $f(x\k)-f(\xstar)$.

In particular, we consider any sequence $\{B\k\}$ that satisfies
\begin{equation}
\label{eq:Bk}
0 \le B\k \leq 2L(\mu/L - \rho)\dist\k, \quad k=0,1,\ldots,
\end{equation}
where $\rho \leq \mu/L$ is a positive constant, which controls the
convergence rate, and $\dist\k$ is a non-negative lower bound on the
distance to optimality, i.e.,
\begin{equation}\label{eq:5}
0 \leq \dist\k \leq f(x\k) - f(\xstar).
\end{equation}

\begin{btheorem}[Strong linear rate under absolute error bounds]
  \label{th:linear-strong}
  Suppose that $\norm{e\k}^2\le B\k$
  where $B\k$ is given by~\eqref{eq:Bk}. Then at each iteration of
  algorithm~\eqref{eq:2} with $\alpha\k\equiv 1/L$,
\begin{equation*}
 f(x\kp1) - f(\xstar) \leq (1-\rho)[f(x\k) - f(\xstar)].
\end{equation*}
\end{btheorem}
\begin{proof}
  Because $\norm{e\k}^2\le B\k$, Lemma~\ref{le:upper-bound}
  and~\eqref{eq:5} imply that
\begin{align*}
\label{eq:xkp1}
f(x\kp1) - f(\xstar) &\leq (1 - \mu/L)[f(x\k) - f(\xstar)] + \frac{1}{2L}B\k\\
& \leq (1 - \mu/L)[f(x\k) - f(\xstar)] + (\mu/L - \rho)\dist\k\\
& \leq (1 - \mu/L)[f(x\k) - f(\xstar)] + (\mu/L - \rho)[f(x\k) - f(\xstar)]\\
& = (1 - \rho)[f(x\k) - f(\xstar)],
\end{align*}
as required.
\end{proof}

As we did with Theorem~\ref{th:linear-rate-weak-expected}, we consider
the case where the approximate gradient is computed by a stochastic
algorithm, and the error can only be bounded in expectation.  Provided
we now have a lower bound $\dist\k$ on the expected sub-optimality,
i.e.,
\[
0 \leq \dist\k \leq \mathbb{E}[f(x\k) - f(\xstar)],
\]
it is possible to show an expected linear convergence rate that
parallels Theorem~\ref{th:linear-strong}. The proof follows that of
Theorem~\ref{th:linear-strong}, where we instead begin by taking the
expectation of both sides of~\eqref{eq:linCon1}.

\begin{btheorem}[Strong expected linear rate under expected error bounds]
  \label{th:linear-strong-expectation}
  Suppose that $\mathbb{E}[\norm{e\k}^2]\le B\k$
  where $B\k$ is given by~\eqref{eq:Bk}. Then at each iteration of
  algorithm~\eqref{eq:2} with $\alpha\k\equiv 1/L$,
  \begin{equation*}
    \mathbb{E}[f(x\kp1) - f(\xstar)] \leq (1-\rho)\mathbb{E}[f(x\k) - f(\xstar)].
  \end{equation*}
\end{btheorem}

In both scenarios, we obtain the fastest convergence rate in the
extreme case where $\rho = \mu/L$.  From~\eqref{eq:Bk}, this means
that $B\k$ must be zero, i.e., the gradient is exact, and we obtain
the classic strong-linear convergence result with error constant
$\sigma=1-\mu/L$ stated in \S\ref{sec:related-work}.  However,
if we take any positive $\rho$ less than $\mu/L$, then we obtain a
slower linear convergence rate but~\eqref{eq:Bk} allows $B\k$ to be
non-zero (as long as $\dist\k > 0$).

The bound~\eqref{eq:Bk} on the error depends on both the conditioning
of the problem (as determined by the parameter $\mu$ and $L$), and on
the lower bound on the distance to the optimal function value
$\dist\k$.  This seems intuitive.  For example, we can allow a larger
error in the gradient calculation the further the current iterate is
from the optimal function value, but a more accurate calculation is
needed to maintain the strong linear convergence rate as the iterates
approach the solution.  Similarly, if the problem is well conditioned
so that the ratio $\mu/L$ is close to 1, a larger error in the
gradient calculation is permitted, but for ill-conditioned problems
where $\mu/L$ is very small, we require a more accurate gradient
calculation.

Note that the analysis of this section holds even if the bounds
$\mu$, $L$, and $\dist\k$ are not the tightest possible. Unsurprisingly,
we obtain the fastest convergence rate when $\mu$ is as large and $L$
is as small as possible, while the largest error in the gradient
calculation is allowed if $\dist\k$ is similarly as large as possible.  Note that
with $\dist\k = 0$, which is a trivial bound that is valid by definition,
we require an exact gradient.

Although it is difficult in general to obtain a bound $\dist\k$ that
satisfies~\eqref{eq:5}, there are at least two possible heuristics
that could be used in practice to approximate this quantity.  The
first heuristic is possible if we have bounds on the Lipschitz and strong-convexity constants
($L$ and $\mu$), as well as the norm of the gradient at $x_k$. In
that case, we use the stationarity of $\xstar$, \eqref{eq:12},
and~\eqref{eq:xstar_fxstar} to deduce that
\[
 \frac{\mu}{2L^2}\norm{\nabla f(x\k)}^2
  \leq \frac{\mu}{2}\norm{x\k - \xstar}^2
  \leq f(x\k)-f(\xstar).
\]
Although we typically do not have access to $\norm{\nabla f(x\k)}$,
the norm of its approximation $\norm{g\k}$ is a reasonable proxy,
which will improve as the magnitude of the gradient residual
decreases.

The second heuristic is based on the assumption that the distance of
consecutive iterates $x_k$ to the solution $x_*$ decreases
monotonically. In that case,
\[
     \norm{x\k - x\kp1}^2
=    \norm{(x\k-\xstar)-(x\kp1-\xstar)}^2
\leq 4\norm{x\k-\xstar}^2.
\]
Coupling this with~\eqref{eq:xstar_fxstar} and premultiplying by
$\mu/8$, we obtain the bound
\[
  \frac\mu8\norm{x\k-x\kp1}^2 \le f(x\k) - f(\xstar).
\]
This option may give a reasonable heuristic in the context of the
increasing sample-size strategy even if the distance to the optimal solution
does increase on some iterations, because as the sample size increases
it becomes less likely that the distance will increase.

\subsection{Relaxing strong convexity}
\label{sec:weakly}

If we remove the assumption of strong convexity, the deterministic
gradient method has a sublinear convergence rate of
$\mathcal{O}(1/k)$~\cite[\S2.1.5]{nesterov2004introductory} while the
stochastic gradient method under standard assumptions has a slower
sublinear convergence rate of
$\mathcal{O}(1/\sqrt{k})$~\cite[\S14.1]{nemirovski1994efficient}.
Using an argument that does not rely on strong convexity, we show that
the $\mathcal{O}(1/k)$ convergence rate of the deterministic gradient
method is preserved for the average of the iterates if the residuals
of the computed gradients are summable.

\begin{btheorem}[Sublinear rate under summable error bounds]
  Suppose that $\sum_{i=0}^\infty \norm{e\k} <
  \infty$. Then at each iteration of algorithm~\eqref{eq:2} with
  $\alpha\k\equiv 1/L$,
\begin{equation*}
 f(\xbar\k) - f(\xstar) = \mathcal{O}(1/k),
\end{equation*}
where $\xbar\k := (1/k)\sum_{i=1}^k x_i$.
\end{btheorem}
\begin{proof}
As in~\eqref{eq:6}, Lipschitz continuity of the gradient implies that
\[
f(x\kp1) \leq f(x\k) - \frac{1}{2L}\norm{\nabla f(x\k)}^2 + \frac{1}{2L}\norm{e\k}^2.
\]
We use convexity of $f$ to bound $f(x\k)$ and obtain
\[
f(x\kp1) \leq f(\xstar) + (x\k-\xstar)^T\nabla f(x\k) - \frac{1}{2L}\norm{\nabla f(x\k)}^2 + \frac{1}{2L}\norm{e\k}^2.
\]
Combine~\eqref{eq:2} and~\eqref{eq:3} to deduce that $-\nabla f(x\k) =
e\k + L(x\kp1 - x\k)$, and move $f(\xstar)$ to the left-hand side to
get, after simplifying,
\begin{align*}
f(x\kp1) - f(\xstar)
& = -(x\k-\xstar)^T(e\k + L[x\kp1 - x\k])\\
&\quad - \frac{L}{2}\norm{x\kp1-x\k}^2 - (x\kp1-x\k)^T e\k\\
& = \frac{L}{2}\norm{x\k-\xstar}^2 -(x\k-\xstar)^T e\k- (x\kp1-x\k)^T e\k\\
&\quad - \frac{L}{2}\norm{x\k-\xstar}^2 - L(x\k-\xstar)^T(x\kp1 - x\k) - \frac{L}{2}\norm{x\kp1-x\k}^2\\
& \leq \frac{L}{2}\left[\norm{x\k-\xstar}^2 - \norm{x\kp1-\xstar}^2\right] + \norm{e\k}\cdot\norm{x\kp1-\xstar}.
\end{align*}
Summing both sides up to iteration $k$ we get
\begin{equation} \label{eq:14}
\sum_{i=1}^k\big[f(x_i) - f(\xstar)\big]
 \leq \frac{L}{2}\norm{x_0-\xstar}^2 - \frac{L}{2}\norm{x\k -
   \xstar}^2 + \sum_{i=1}^k\norm{e_{i-1}}\cdot\norm{x_i-\xstar}.
\end{equation}
We will first use~\eqref{eq:14} to show that the sequence
$\{\norm{x\k-\xstar}\}$ is bounded, and then use it to obtain the final result.

Because the left-hand side of~\eqref{eq:14} is nonnegative,
\begin{equation*}
\frac{L}{2}\norm{x\k-\xstar}^2 \le \frac{L}{2}\norm{x_0-\xstar}^2 + \sum_{i=1}^k\norm{e_{i-1}}\cdot\norm{x_i-\xstar}.
\end{equation*}
We now prove that the sequence
$\{\norm{x\k-\xstar}\}$ is bounded by showing that the
auxiliary sequence $\{d\k\}$, with $d\k := \max\{1,\norm{x\k-\xstar}\}$, is
bounded. Because $d\k\ge1$,
\begin{align*}
  d\k \le d\k^2
    &\le \max\left\{\, 1,\, \norm{x_0-\xstar}^2
                     +
              \frac{2}{L}\sum_{i=1}^k\norm{e_{i-1}}\cdot\norm{x_i-\xstar}
              \, \right\}
\\ &\le d_0^2 + \frac{2}{L}\sum_{i=1}^k\norm{e_{i-1}}\cdot\norm{x_i-\xstar}
\\ &\le d_0^2 + C_1\sum_{i=1}^k\norm{e_{i-1}}\cdot d_i,
\end{align*}
where $C_1:=2/L$.
Because the sequence $\{\norm{e\k}\}$ is summable, it holds that
$\norm{e\k}\to0$, and thus that there exists some $N$ large enough
 that $C_1\norm{e_{k-1}}<1$ for all $k\ge N$. Partition the sum at $N$:
\[
  d\k \le C_0 + C_1\sum_{i=N}^k\norm{e_{i-1}}\cdot d_i
  \text{for all} k>N,
\]
where $C_0:=d_0^2+C_1\sum_{i=1}^{N-1}\norm{e_{i-1}}\cdot
d_i$. Rearrange terms to get
\[
 (1-C_1\norm{e\km1})d\k \le C_0 + C_1\sum_{i=N}^{k-1}\norm{e_{i-1}}\cdot d_i
 \text{for all} k > N,
\]
and because $C_1\norm{e_{k-1}}<1$ for all $k>N$,
\begin{equation*} 
 d\k \le \frac{1}{1-C_1\norm{e_{k-1}}}
 \left[
   C_0 + C_1\sum_{i=N}^{k-1}\norm{e_{i-1}}d_i
 \right]
 \text{for all}
 k > N.
\end{equation*}
If we apply this bound recursively to $d_{k-1}$ we obtain
\begin{align*}
 d\k & \le \frac{1}{1-C_1\norm{e_{k-1}}}
 \hspace{-1pt}
 \left[
   C_0 + C_1\sum_{i=N}^{k-2}\norm{e_{i-1}}d_i + \frac{C_1\norm{e_{k-2}}}{1-C_1\norm{e_{k-2}}}
   \hspace{-1pt}\left(C_0 + C_1\sum_{i=N}^{k-2}\norm{e_{i-1}}d_i\right)
 \right]\\
 & = \frac{1}{(1-C_1\norm{e_{k-1}})(1-C_1\norm{e_{k-2}})}\left[C_0 + C_1\sum_{i=N}^{k-2}\norm{e_{i-1}}d_i\right]
 \text{for all}
 k > N,
\end{align*}
and if we apply it recursively from $k$ down to $N$ we obtain
\begin{equation} \label{eq:28}
  d\k \le \frac{2C_0}{\prod_{i=N}^k     \left(1-C_1\norm{e_{i-1}}\right)}
      \le \frac{2C_0}{\prod_{i=N}^\infty\left(1-C_1\norm{e_{i-1}}\right)}.
\end{equation}
To see that the right-hand side term is bounded, take
logarithms of both sides to get
\[
  \log d\k \le \log(2C_0) - \sum_{i=N}^\infty\log(1-C_1\norm{e_{i-1}}).
\]
We now use the limit-comparison test, which asserts that for nonnegative
sequences $\{a_i\}$ and $\{b_i\}$, if $0< \lim_{i\to\infty} a_i/b_i <
\infty$, then $\sum_i b_i<\infty$ implies $\sum_i a_i<\infty$. We thus
compare the sequence $\{-\log(1-C_1\norm{e_{i-1}})\}$ to the summable
sequence $\{\norm{e_{i-1}}\}$ using L'H\^opital's rule to get
\[
  \lim_{i\to\infty} \frac{-\log(1-C_1\norm{e_{i-1}})}{\norm{e_{i-1}}}
  = \lim_{i\to\infty} \frac{\frac{C_1}{1-C_1\norm{e_{i-1}}}}{1}
  = C_1.
\]
Thus, the sequence $\{-\log(1-C_1\norm{e_{i-1}})\}$ is summable, which
implies, via~\eqref{eq:28} and the definition of $d\k$, that the sequence $\{\norm{x\k-\xstar}\}$
is bounded.

We now use convexity of $f$ and~\eqref{eq:14}  to bound the function value
of the average iterate:
\begin{align*}
f\left(\frac{1}{k}\sum_{i=1}^kx_i\right)-f(\xstar) 
   & \leq \frac{1}{k}\sum_{i=1}^k \big[f(x_i) - f(\xstar)\big]
\\ & \leq \frac{L}{2k}\norm{x_0 - x_*}^2
  + \frac{1}{k}\sum_{i=1}^k\norm{e_{i-1}}\cdot\norm{x_i-\xstar}.
\end{align*}
Because the sequence $\{\norm{e_k}\}$ is summable (by assumption), and
the sequence $\{\norm{x_k-\xstar}\}$ is bounded, this last inequality
implies the conclusion of the theorem.
\end{proof}

Note that the convergence rate also holds for the iterate that
achieves the lowest function value but, unlike the deterministic case,
this is not guaranteed to be the last iteration.

\section{Application to sample-average gradients}
\label{sec:batch-error}

Incremental-gradient methods for~\eqref{eq:1} are based on the
iteration scheme~\eqref{eq:2} with the gradient approximation
\begin{equation*}
  g\k \defd \frac1{\abs{\batch\k}}\sum_{i \in\batch\k} \nabla f_i(x\k),
\end{equation*}
where the set $\batch_k\subseteq\{1,\dots,M\}$ represents a sample of
the measurements that constitute the full data set.
Typically, $\batch\k$ contains a single element that is chosen in
either a cyclic fashion or sampled at random, and as discussed
in~\S\ref{sec:related-work}, the convergence rate of the method is
sublinear.  As the sample size increases, however, the error in the
sampled gradient $g\k$ decreases, and so the sample size can be used
to implicitly control the error in the gradient.  We use the results
of~\S\ref{sec:convergence} to develop an increasing sample-size
strategy that improves on this sublinear rate.

Let $\batchcomp\k$ denote the complement of $\batch\k$, so that
$\batch\k\cup\batchcomp\k= \{1,\dots,M\}$. Then the gradient residual,
defined by~\eqref{eq:3}, satisfies
\begin{equation}\label{eq:23}
e\k = \frac{M-\abs{\batch\k}}{M\abs{\batch\k}}\sum_{i \in \batch\k} \nabla f_i(x\k)
    - \frac{1}{M}\sum_{i\in\batchcomp\k} \nabla f_i(x\k).
\end{equation}
The first term is a re-weighting of the gradient approximation and the 
second term is the portion of $\nabla f(x\k)$ that is not
sampled. Assumption~\eqref{eq:incGradBnd} allows us to bound the norm
of this residual in terms of the full gradient, and thus leverage the
results of~\S\ref{sec:convergence}.

\subsection{Deterministic sampling}
\label{sec:deterministic_error_bounds}

It follows from assumption~\eqref{eq:incGradBnd} that the gradient
residual in~\eqref{eq:23} satisifies
\begin{align*}
  \norm{e\k}^2
  &=
  \left\|
    \left(\frac{M-\abs{\batch\k}}{M\abs{\batch\k}}\right)
    \sum_{i\in\batch\k}\nabla f_i(x\k)
    - \frac1M\sum_{i\in\batchcomp\k}\nabla f_i(x\k)
  \right\|^2
  \\ &\le
  \left[
    \left(
      \frac{M-\abs{\batch\k}}{M\abs{\batch\k}}
    \right)
      \Big\|\sum_{i\in\batch\k}\nabla f_i(x_k)\Big\|
      +
      \frac{1}{M}\Big\|\sum_{i\in\batchcomp\k}\nabla f_i(x\k)\Big\|
 \right]^2
  \\ &\le
  \left[
    \left(
      \frac{M-\abs{\batch\k}}{M\abs{\batch\k}}
    \right)
      \sum_{i\in\batch\k}\norm{\nabla f_i(x_k)}
      +
      \frac{1}{M}\sum_{i\in\batchcomp\k}\norm{\nabla f_i(x\k)}
 \right]^2
  \\ &\le
  4\left(
    \frac{M-\abs{\batch\k}}{M}
  \right)^2
  (\beta_1+\beta_2\norm{\nabla f(x\k)}^2),
\end{align*}
where the triangle-inequality is applied repeatedly, and the resulting
terms simplified. Next, in the same way that~\eqref{eq:grdBnd} was
derived, we use~\eqref{eq:lip} to derive the upper bound
\[
 \norm{\nabla f(x\k)}^2 \le 2L[f(x\k)-f(\xstar)].
\]
Thus the bound on $\norm{e\k}^2$ can be expressed in terms of
the sample-size ratio and the distance to optimality:
\begin{equation}\label{eq:7}
 \norm{e\k}^2 \le 4\left[\frac{M-\abs{\batch\k}}{M}\right]^2
 \big(\beta_1+2\beta_2L[f(x\k) - f(\xstar)]\big).
\end{equation}

The following result parallels Theorem~\ref{th:linear-rate-weak}, and
asserts that a linearly increasing sample size is sufficient to induce a
weak linear convergence rate of the algorithm.
\begin{btheorem}[Weak linear rate with deterministic sampling]
  \label{th:incremental-weak}
  Suppose that \eqref{eq:incGradBnd} holds, and that the sample size
  $\abs{\batch\k}$ increases geometrically towards $M$, i.e.,
  \begin{equation*}
    \label{eq:9}
    \frac{M - \abs{\batch\k}}{M} = \Oscr(\gamma^{k/2})
  \end{equation*}
  for some $\gamma<1$.   
  Then at each iteration of
  algorithm~\eqref{eq:2} with $\alpha\k\equiv 1/L$, for any
  $\epsilon>0$ we have
  \begin{equation*}
  f(x\k) - f(\xstar) = [f(x_0)-f(x_*)]\Oscr([1-\mu/L+\epsilon]^k) + \Oscr(\sigma^k)
  \end{equation*}
  where $\sigma = \max\{\gamma,\, 1-\mu/L\}+\epsilon$.
\end{btheorem}
\begin{proof}
  Let $\rho\k=\big[\frac{M-\abs{\batch\k}}{M}\big]^2$. Using (3.2) and
  Lemma~2.1, we obtain the bound
  \begin{equation*}
    \begin{aligned}
      f(x\kp1) - f(\xstar)
        &\le (1-\mu/L)[f(x\k)-f(\xstar)] +
             \frac{2\rho\k}{L}(\beta_1 + 2\beta_2L[f(x\k)-f(\xstar)])
      \\&=   (1-\mu/L+4\beta_2\rho\k)[f(x\k)-f(\xstar)]
             + \frac{2\beta_1}{L}\rho\k
      \\&= \omega\k[f(x\k) - f(\xstar)] + \frac{2\beta_1}{L}\rho\k,
    \end{aligned}
  \end{equation*}
  where $\omega\k\defd 1-\mu/L+4\beta_2\rho\k$. Apply this
  recursively and use $\rho_k = \Oscr(\gamma^k)$ to obtain
  \begin{align*}
  f(x\k) - f(\xstar)
  &\le [f(x_0)-f(x_*)]\prod_{i=0}^{k-1}\omega_i
      + \sum_{i=0}^{k-1}\Oscr(\gamma^i\prod_{j=i+1}^{k-1}\omega_j).
  \end{align*}
  Take $\delta\k \defd \max\{\gamma,\,\omega\k\}$. Because
  $\gamma^i = \prod_{j=1}^i\gamma \leq \prod_{j=1}^i\delta_j$
for
  all $i$,
  \begin{align*}
  f(x\k) - f(\xstar)
  &\le [f(x_0)-f(x_*)]\prod_{i=0}^{k-1}\omega_i +
    \sum_{i=0}^{k-1}\Oscr\Big(\prod_{j=1}^{k-1}\delta_j\Big)
  \\
  &= [f(x_0)-f(x_*)]\prod_{i=0}^{k-1}\omega_i + \Oscr\Big(k\prod_{j=1}^{k-1}\delta_j\Big).
  \end{align*}
  Because $\rho\k\to0$, it follows that $\omega\k\to1-\mu/L$ and thus
  $\prod_{i=0}^k\omega_i=\Oscr([1-\mu/L+\epsilon]^k)$ for any
  $\epsilon>0$, which bounds the first term in the right-hand side
  above. Furthermore, we now use the fact that
  $\delta\k\searrow\deltabar\defd\max\{\gamma,\,1-\mu/L\}$ to show
  that the second term is in $\Oscr([\deltabar+\epsilon]^k)$ for any
  $\epsilon>0$.
  In particular, choose $N$ large enough that
   $(\delta_k+\delta_k/k) \leq \deltabar + \epsilon$ for all $k \geq N$
   and choose the constant $\xi$ so that
   \[
   k\prod_{j=1}^{k-1}\delta_j \leq \xi(\deltabar + \epsilon)^k,
   \]
   for all $k < N$. Then by induction, $k\prod_{j=1}^{k-1}\delta_j = \Oscr(\sigma^k)$ for all $k$ because
   \begin{align*}
   (k+1)\prod_{j=1}^k\delta_j & = (\delta_k+\delta_k/k)k\prod_{j=1}^{k-1}\delta_j \\
 & \leq (\delta_k+\delta_k/k)\xi(\deltabar + \epsilon)^k \\
 & \leq  (\deltabar + \epsilon)\xi(\deltabar + \epsilon)^k\\
 & = \xi(\deltabar + \epsilon)^{k+1}.
     \end{align*}
    \end{proof}

An interesting difference with Theorem~\ref{th:linear-rate-weak},
which considers a generic error in the gradient, is that the error in
the objective function decreases at twice the rate that the sample
size increases.

It is possible to get a strong linear rate of convergence by
increasing the sample size in a more controlled
way. Theorem~\ref{th:linear-strong} guides the choice of the sample
size, and gives the following corollary. The proof follows by simply
ensuring that the right-hand side of~\eqref{eq:7} is bounded as
required by Theorem~\ref{th:linear-strong}.
\begin{bcorollary}[Strong linear rate with deterministic sampling]
  \label{th:incremental-strong}
  Suppose that \eqref{eq:incGradBnd} holds, and that the sample size
  $\abs{\batch\k}$ is increased so that at each iteration
  $k=0,1,\ldots$,
  \begin{equation*} \label{eq:11}
    4\left(\frac{M-\abs{\batch\k}}{M}\right)^2
    \big(\beta_1+2\beta_2L[f(x\k) - f(\xstar)]\big)
    \le
    L(\mu/L-\rho)[f(x\k)-f(\xstar)]
  \end{equation*}
  for some positive $\rho\le\mu/L$. Then at each iteration of
  algorithm~\eqref{eq:2} with $\alpha\k\equiv 1/L$,
  \begin{equation*}
  f(x\kp1) - f(\xstar) \le (1-\rho)[f(x\k) - f(\xstar)].
  \end{equation*}
\end{bcorollary}

We see that if the individual functions $f_i$ are very similar (so
that $\beta_1$ and $\beta_2$ are small) then we can choose a fairly
small sample size.  In contrast, if the $f_i$ are very dissimilar then
we must use a larger sample size.  

\subsection{Stochastic sampling}

Theorem~\ref{th:incremental-weak} and Corollary~\ref{th:incremental-strong} are based on the deterministic
bound~\eqref{eq:7} on the gradient error, and hold irrespective of
the manner in which the elements of the samples $\batch\k$ are chosen,
e.g., the samples do not need to be chosen cyclically or sampled
uniformly.  We can obtain a strictly tighter bound in expectation, however, if we choose the
sample by uniform sampling \emph{without replacement}.
In particular, by suitably modifying the derivation of a well-known
result in statistical sampling, we can obtain a bound in terms of a quantity
related to the sample variance of the gradients:
\[
S\k = \frac{1}{M-1}\sum_{i=1}^M||\nabla f_i(x\k) - \nabla f(x\k)||^2.
\]
Using an argument similar to~\cite[\S2.7]{lohr1999sampling},
uniform sampling \emph{without replacement} yields
\[
  \mathbb{E}[\norm{e\k}^2] =
  \left(
    \frac{M-\abs{\batch\k}}{M}
  \right)\frac{S\k}{\abs{\batch\k}}.
\]
By using~\eqref{eq:incGradBnd}, we obtain the bound
\begin{equation*}
 \eval[\norm{e\k}^2]
 \leq
   \left[
     \frac{M-\abs{\batch\k}}{M}
     \cdot \frac1{\abs{\batch\k}}
   \right]
   \frac{M}{M-1}(\beta_1 +
 2(\beta_2-1)L[f(x\k)-f(\xstar)]).
\end{equation*}
An expected convergence result parallel to
Theorem~\ref{th:incremental-weak} can then be obtained with only minor
changes to the proof.
\begin{btheorem}[Weak linear rate with stochastic sampling]
  \label{th:randomized-incremental-weak}
  Suppose that \eqref{eq:incGradBnd} holds, and that
  \begin{equation*}
    \frac{M - \abs{\batch\k}}{M}\cdot\frac{1}{\abs{\batch\k}} = \Oscr(\gamma^k)
  \end{equation*}
  for some $\gamma<1$.   
  Then at each iteration of
  algorithm~\eqref{eq:2} with $\alpha\k\equiv 1/L$, for any
  $\epsilon>0$ we have
  \begin{equation*}
  \mathbb{E}[f(x\k) - f(\xstar)] = [f(x_0)-f(x_*)]\Oscr([1-\mu/L+\epsilon]^k) + \Oscr(\sigma^k)
  \end{equation*}
  where $\sigma = \max\{\gamma,\, 1-\mu/L\}+\epsilon$.
\end{btheorem}

\begin{figure}[t]
  \centering\small
  \begin{tabular}{@{}cc@{}}
  \includegraphics[width=.47\textwidth]{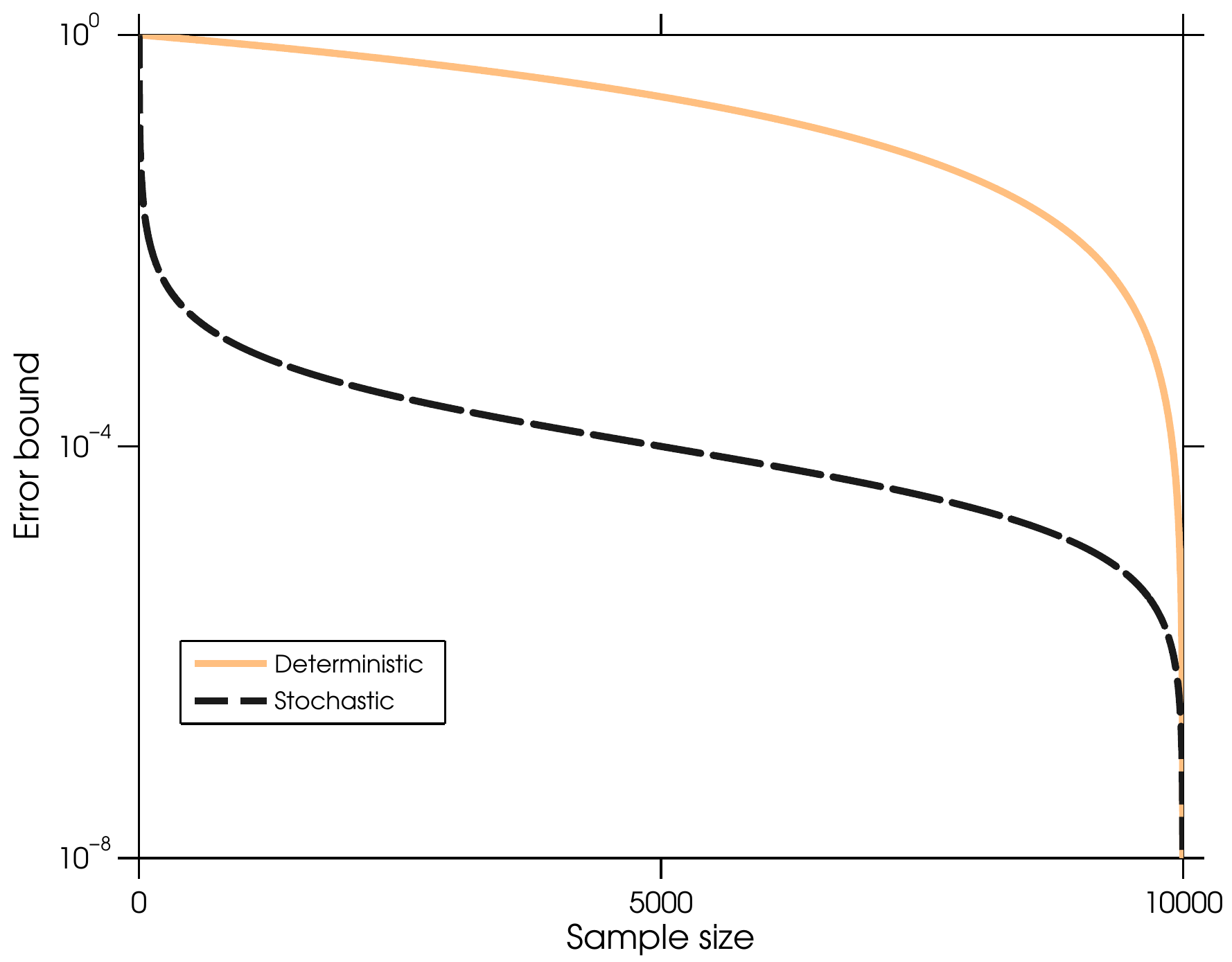}
 &\includegraphics[width=.47\textwidth]{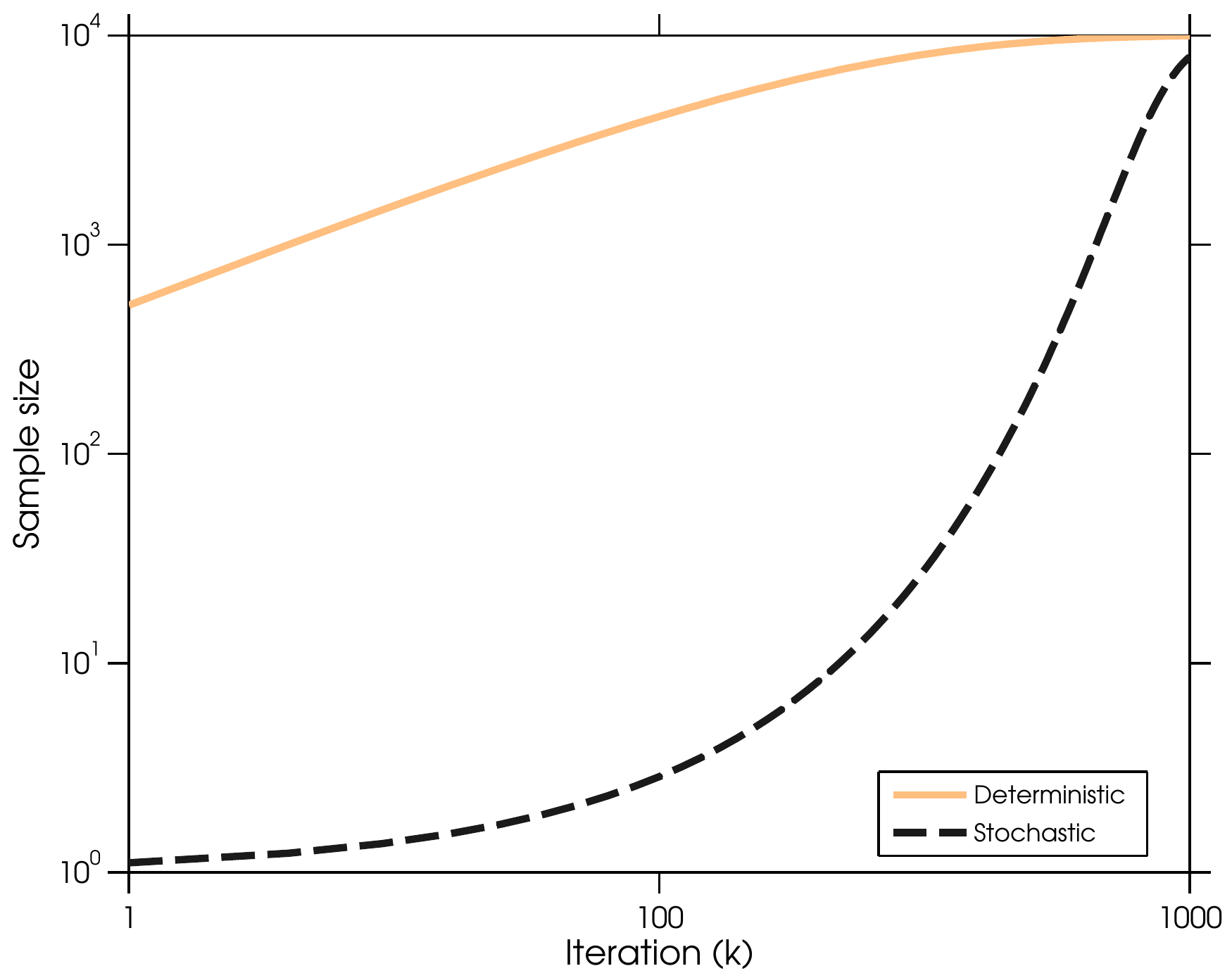}
\\(a) & (b)
  \end{tabular}
  \caption{(a) Bounding factors in the error of the sampled
    gradient. The error bound for stochastic sampling is uniformly
    better than for deterministic sampling. (b) Minimum sample-size
    schedule required to achieve a linear rate with error constant
    0.9. In both figures, $M=10^4$.}
  \label{fig:error-bnd}
\end{figure}

Note that the right-hand side of the bound on
$\mathbb{E}[\norm{e\k}^2]$ is uniformly better than the bound shown
in~\eqref{eq:7}. Importantly, it initially decreases to zero at a
faster rate as the size of the sample
increases. Figure~\ref{fig:error-bnd}(a) illustrates the difference in
the sample-size requirements between
Theorems~\ref{th:incremental-weak}
and~\ref{th:randomized-incremental-weak}, i.e.,
\[
\left[\frac{M-\abs{\batch\k}}{M} \right]^2
\quad \hbox{(deterministic)}
\textt{vs.}
\left[\frac{M-\abs{\batch\k}}{M}\cdot\frac{1}{\abs{\batch\k}}\right]
\quad \hbox{(stochastic)},
\]
as the sample size $\abs{\batch\k}\to
M:=10^4$. Figure~\ref{fig:error-bnd}(b) illustrates the sample-size
schedule needed to realize a linear convergence rate in the
deterministic and stochastic cases.

\section{Practical implementation}
\label{sec:practicalImplementation}

The analysis so far has focused on the approximate gradient descent
iterations given by \eqref{eq:2}--\eqref{eq:3}. In practice, it is
useful to make two modifications to this basic algorithm. The first,
described by \eqref{eq:25}--\eqref{eq:27}, is based on scaling the
search directions to account for curvature information in $f$. The
second allows for a varying step size. For this section, we define the
sampled objective function and its gradient by
\begin{equation} \label{eq:29}
  \fbar\k(x) = \frac{1}{\abs{\batch\k}}\sum_{i\in\batch\k}f_i(x)
  \text{and} 
  g\k(x) = \frac{1}{\abs{\batch\k}}\sum_{i\in\batch\k}\nabla f_i(x).
\end{equation}

\subsection{Scaled Direction}
\label{sec:scaledDirection}

Our implementation attempts to gather useful curvature information on
the function $f$ by maintaining a quasi-Newton approximation to the
Hessian $H\k$. The search directions $d\k$ are then taken as the solution of
the system~\eqref{eq:27}.
The approximate Hessian is maintained via the recursive application of the
update
\[
H\kp1 = U(H\k, y\k, s\k),
\]
where $U$ represents an update formula on the $k$th iteration, while
\[
 s\k \defd x\kp1 - x\k \text{and} y\k \defd g\k(x\kp1) - g\k(x\k)
\]
measure the change in $x$ and the sampled gradient. In our experiments
(see \S\ref{sec:experiments}) we use a limited-memory quasi-Newton
update, which maintains a history of the previous $\ell=10$ pairs
$(s\k,y\k)$, and recursively applies the formula $H_{i+1} =
U(H_i,s_i,y_i)$ for $i=k-\ell,\ldots,k$. Nocedal and
Wright~\cite[\S7.2]{NoceWrig:2006} 
describe the recursive procedure for the limited-memory BFGS update
that we use. Updates are skipped if necessary in order to ensure that
the approximation $H\k$ remains positive definite and has a bounded
condition number. To ensure that $d\k$ is well-scaled, we use the
Shanno and Phua~\cite{ShannoPhua:1978} scaling of the initial Hessian
approximation on each iteration as described by Nocedal and
Wright~\cite[p.~178]{NoceWrig:2006}.

Another approach, not considered here, is to base the quasi-Newton
Hessian on an independent set of sampled gradients~\cite{BCNN:2011}.




\subsection{Varying stepsize}
\label{sec:stepsize}

A weakness of our convergence analysis is the requirement for a fixed
steplength $\alpha\k\equiv L$, in part because the Lipschitz constant
is not usually known, and in part because a dynamic steplength is
typically more effective in practice. But a linesearch
procedure that ensures a sufficient decrease condition in the true
function $f$ runs contrary to a sampling scheme specifically
designed to avoid expensive evaluations with $f$.

In our implementation we attempt to strike a balance between a
rigorous linesearch and none at all by enforcing an Armijo-type
descent condition on the sampled objective $\fbar\k$. In
particular, 
we
use a linesearch procedure to select a steplength $\alpha\k$ that
satisfies
\begin{equation}  \label{eq:13}
  \fbar\k(x\k + \alpha d\k ) < \fbar\k(x\k) + \eta\alpha g\k(x\k)\T d\k,
\end{equation}
where $d\k$ is the current search direction and $\eta\in(0,1)$.
While in deterministic quasi-Newton methods we typically first test
whether $\alpha = 1$ satisfies this condition, in our implementation
we set our initial trial step length to $\alpha = |\batch_{k-1}|/|\batch_k|$.

In general, $\batch\k$ represents only a fraction of all observations,
and so the above procedure may not even yield a decrease in the true
objective at each iteration. But because we steadily increase the
sample size, the procedure just described eventually reduces to a
conventional linesearch based on the true objective function
$f$. Hence, the method may initially be nonmonotic, but is guaranteed
to eventually be monotonic. Coupled with the choice of search
direction described in \S\ref{sec:scaledDirection}, the overall
algorithm reduces to a conventional linesearch method with a
quasi-Newton Hessian approximation, and inherits the global and local
convergence guarantees of that method.

\section{Numerical experiments}
\label{sec:experiments}

This section summarizes a series of numerical experiments in which we
apply our incremental-gradient method with a growing sample size to a
series of data-fitting applications. Table~\ref{tab:test-problems}
summarizes the test problems.

The first four experiments are data-fitting applications of logistic
regression of varying complexity: binary, multinomial,
chain-structured conditional random fields (CRFs), and general
CRFs. These logistic-regression applications follow a standard
pattern. We first model the probability of an outcome $b_i$ by some
log-concave function $p(b_i \mid a_i, x)$, where $a_i$ is data; the goal is
to choose the parameters $x$ so that the likelihood function
$$
\Lscr(x) = \prod_{i=1}^M p(b_i \mid a_i, x)
$$
is maximized. 
We then approximate $x$ by minimizing the
2-norm-regularized negative log-likelihood function
\begin{equation}\label{eq:18}
  f(x) = \sum_{i=1}^M f_i(x) + \half\lambda\norm{x}^2
  \text{with}
  f_i(x) = -\log p(b_i \mid a_i, x)
\end{equation}
for some positive regularization parameter $\lambda$.  These
objectives are all strongly convex (with $\mu \geq \lambda$) and satisfy our assumptions in
\S\ref{sec:assumptions}. For some functions $p$ the 
Lipschitz constant of $\nabla f$, or upper bounds on it, are available.  We discuss
the particular case of the binary model in \S\ref{sec:binary-logist-regression}.

The last experiment is a more general application of nonlinear
least-squares to seismic inversion. This last data-fitting application
does not satisfy our central convexity assumption, but nonetheless
illustrates the practical relevance of our approach on difficult
problems.

\begin{table}[t]
  \centering \small
  \label{tab:test-problems}
  \caption{The test problems.}
  \begin{tabular}{lrrl}
    \toprule
    Problem & $M$ & $n$ & Description
 \\ \midrule
    binary logistic regression & 92,189 & 823,470 & spam
    identification
    (\S\ref{sec:binary-logist-regression})
 \\ multinormal logistic regression & 70,000 & 785 & digit
    identification (\S\ref{sec:mult-logist-regression})
 \\ chain-structured CRF & 8,936 & 1,643,004 & noun-phrase chunking (\S\ref{sec:crf})
 \\ general CRF & 50 & 4 & image denoising (\S\ref{sec:g-crf})
 \\ nonlinear least squares & 101 & 10,201 & seismic inversion (\S\ref{sec:nnls})
 \\ \bottomrule
  \end{tabular}
\end{table}

Our numerical experiments compare the following three
methods:

\paragraph{Deterministic} A conventional quasi-Newton linesearch
method that uses the true function $f$ and gradient $\nabla f$ at
every iteration. The method is based on a limited-memory BFGS Hessian
approximation and a linesearch based on Hermite cubic-polynomial
interpolation and the strong Wolfe conditions.  Several comparison
studies indicate that these type of limited-memory quasi-Newton
methods are among the most efficient deterministic methods available
for solving large-scale logistic regression and conditional random
field
problems~\cite{malouf2002comparison,wallach2002efficient,minka2003comparison,sha2003shallow}.

\paragraph{Stochastic} An incremental-gradient method based on the
iteration~\eqref{eq:2} with $g\k = \nabla f_i(x\k)$ and a constant $\alpha\k$,
where  the index
$i\in\{1,\ldots,M\}$ is randomly selected. This corresponds to a
constant sample size of one, and this simple method has proved
competitive with more advanced deterministic methods like the one
above for estimation in CRF models~\cite{vishwanathan2006accelerated}.

\paragraph{Hybrid} This is the proposed method described in
\S\ref{sec:practicalImplementation}, which uses search directions
computed from \eqref{eq:27} and a linesearch based on satisfying
condition \eqref{eq:13}. As with the deterministic method described
above, the Hessian approximations are based on limited-memory BFGS and
the linesearch uses polynomial interpolation.  The objective and
gradient approximations are based on~\eqref{eq:29}, where
$\batch\k\subseteq\{1,\ldots,M\}$, and the number of elements in the
current sample is initially 1, and grows linearly as per the formula
$$
\abs{\batch\kp1} = \lceil\min\{1.1\cdot\abs{\batch\k} + 1,\, M\}\rceil.
$$
The hybrid nature of this approach should now be clear: the very
first iteration is similar to the stochastic method; when the sample
size grows to include all observations, the algorithm morphs into the
deterministic method.

All experiments are carried out using Matlab R2010b on a 64-bit Athlon
machine. Two plots are shown for each experiment. The first shows
the progress of the objective value against the index $p =
\frac1M\sum_{k=0}^p\abs{\batch\k}$ for $p=1,2,\ldots$, which measures
the effective number of passes through the entire data set. The second
plot shows the cumulative number of $f_i$ functions evaluated, i.e.,
$\sum_{i=0}^k\abs{\batch_i}$, against the iterations $k$.

\subsection{Binary logistic regression} \label{sec:binary-logist-regression}

Logistic regression models~\cite[\S4.3.2]{bishop2006pattern} are used
in an enormous number of applications for the problem of \emph{binary
  classification}.  We are given data with $M$ examples of
input-output pairs $(a_i,b_i)$, where $a_i \in \Real^n$ is a vector of
$n$ features, and $b_i \in \{-1,1\}$ is a corresponding binary
outcome.  The goal is to build a linear classifier that, given the
features $a_i$ and a vector of parameters $x$, the sign of the
inner-product $a_i\T x$ gives $b_i$. The logistic model gives the
probability that $b_i$ takes the value $1$:
\[
  p_1(b_i = 1 \mid a_i,x) = \frac{\exp(a_i\T x)}{\exp(a_i^Tx) + 1}
                          = \frac{1}{1 + \exp(-a_i\T x)}.
\]
(Typically a bias variable is added, but we equivalently assume that
the first element of $x$ is set to one.)  Thus, the probability
that $b_i$ takes the value -1 is $[1-p(b_i=1 \mid a_i, x)]$.  We can
write these two cases compactly as
\[
  p_1(b_i \mid a_i,x) = \frac{1}{1 + \exp(-b_ia_i^Tx)}.
\]
This is the probability function $p$ used in~\eqref{eq:18}.

\begin{figure}[t]
  \centering
  \includegraphics[width=.48\textwidth]{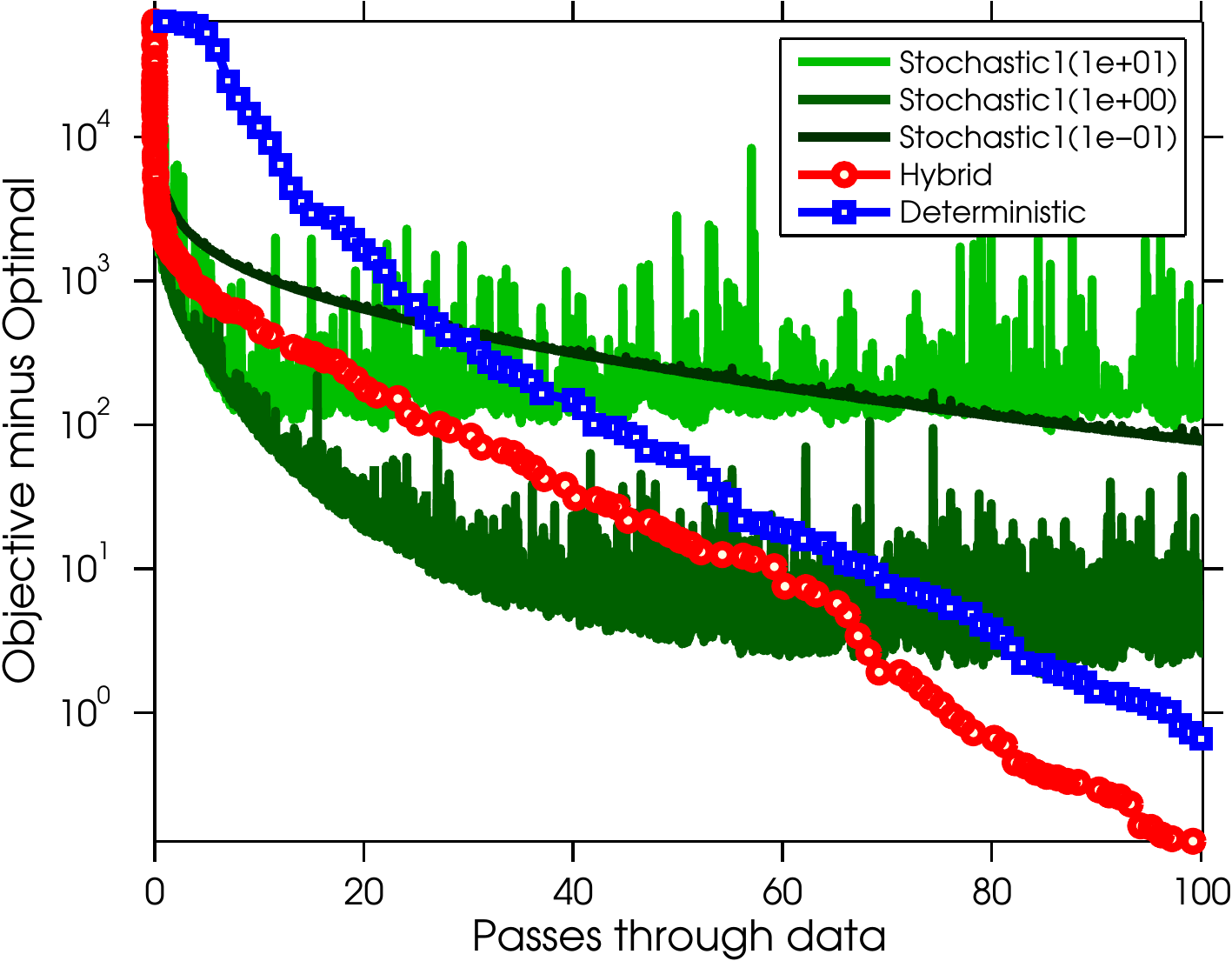}
  \hfill
  \includegraphics[width=.48\textwidth]{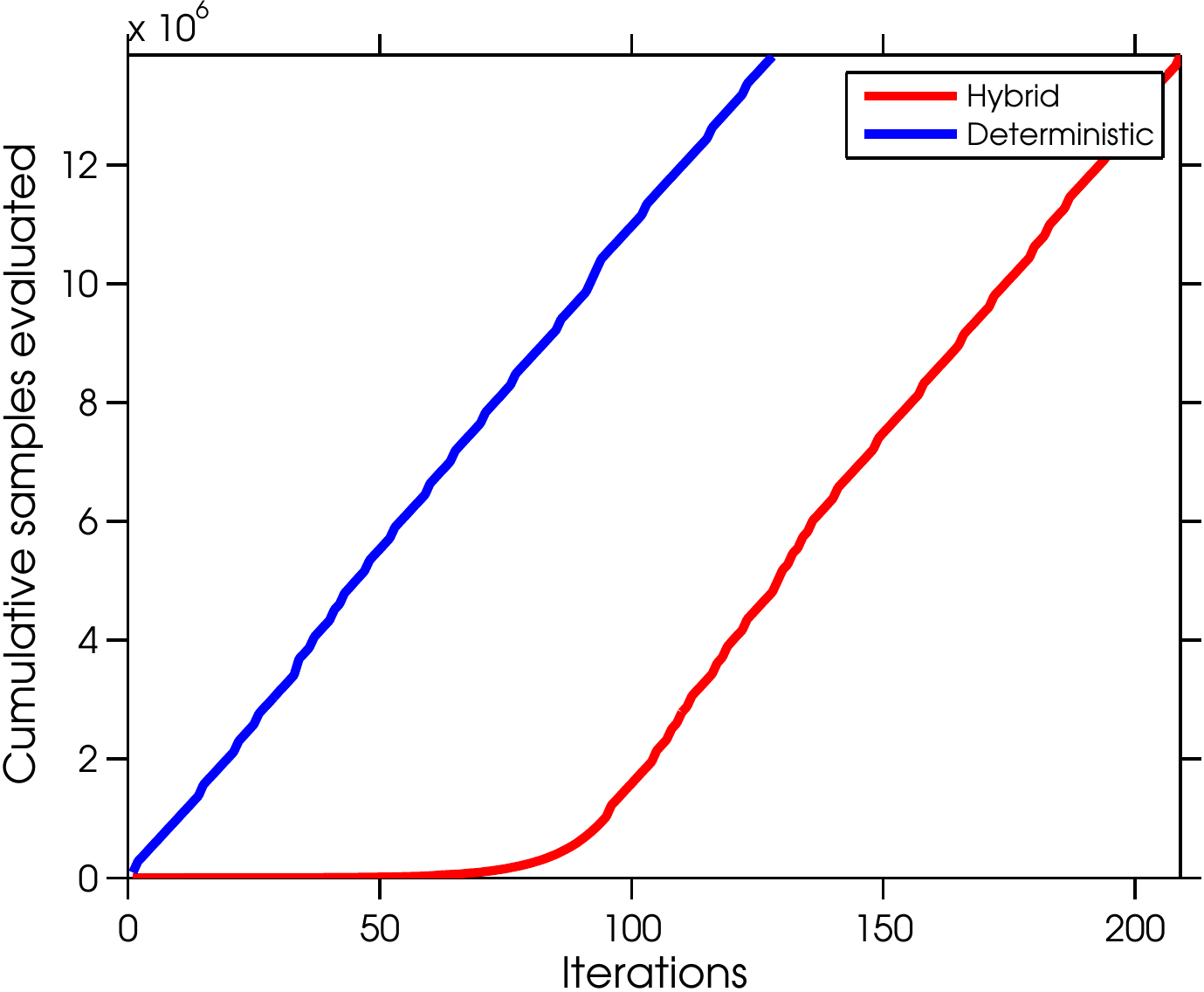}
  \caption{Binary logistic regression experiments for different
    optimization strategies for spam classification. The stochastic
    method is run with 3 different fixed steplengths.}
  \label{fig:logreg}
\end{figure}

The dominant cost in computing $f$ and its gradient is the cost of
forming the matrix-vector products $Ax$ and $A\T y$ (for some $y$), where the $M$
rows of the matrix $A$ are formed from the vectors $a_i$.
The Hessian is $\nabla^2 f(x)=A\T D A$, where
$D$ is a diagonal with elements $p_1(b_i\mid
a_i,x)\cdot[1-p_1(b_i\mid a_i,x)]$, which lie in the range $(0,0.25]$.
The (nonnegative) eigenvalues of the Hessian are thus bounded above by
$0.25\norm{A}^2$, and are strictly positive if $A$ has full rank.
Combined with 2-norm regularization, the resulting function $f$
satisfies the assumptions of~\S\ref{sec:assumptions}.


Our experiments for binary logistic regression are based on the TREC
2005 data set, which contains 823,470 binary variables describing the
presence of word tokens in 92,189 email messages involved in the legal
investigations of the Enron corporation~\cite{cormack2005spam}.  The
target variable indicates whether the email was spam or not.  The
data set was prepared by Carbonetto~\cite[\S2.6.5]{Carbonetto09}, and
we set the regularization parameter $\lambda=0.01$.

The results of this experiment are plotted in Figure~\ref{fig:logreg},
where we define the optimal value as the best value found across the
methods after 150 effective passes through the data, and where for the
stochastic method we plot the three step sizes (among powers of 10)
that gave the best performance over the allotted iterations.  Of
course, in practice we will not know what step size optimizes the
performance of the stochastic method, and the lack of sensitivity of
the result to the initial step size is an advantage of the
deterministic and hybrid methods.  In the plot we see that, like the
stochastic methods, the hybrid method makes rapid initial
progress. Unlike the stochastic method, however, the hybrid method
continues to make steady progress similar to the deterministic method.
This behavior agrees with the theory.

\subsection{Multinomial logistic regression}
\label{sec:mult-logist-regression}

Multinomial logistic regression relaxes the binary requirement, and
allows each outcome $b_i$ to take any value from a set of classes
$\Cscr$~\cite[\S4.3.4]{bishop2006pattern}. In this model there is a
separate parameter vector $x_j$ for each class $j\in\Cscr$. We model
the probability that $b_i$ is assigned a particular class $j$ as
\begin{equation}\label{eq:19}
  p_2(b_i = j \mid a_i, \{x_j\}_{j\in\Cscr}) =
  \frac{\exp(x\j\T a_i\drop)}{\sum_{j'\in\Cscr}\exp(x_{j'}^T a_i\drop)}.
\end{equation}
This model is equivalent to binary logistic regression in the special
case where the parameters $x_j$ of one class are fixed at zero and there are only two classes, i.e., $\Cscr=\{-1,1\}$.
As with binary logistic regression, the function $p_2$ is log-concave and $-\log p_2$ has
a Hessian whose eigenvalues are bounded above.
Hence, the resulting function $f$ in~\eqref{eq:18} satisfies the
assumptions of \S\ref{sec:assumptions}.

\begin{figure}[t]
  \centering
  \includegraphics[width=.48\textwidth]{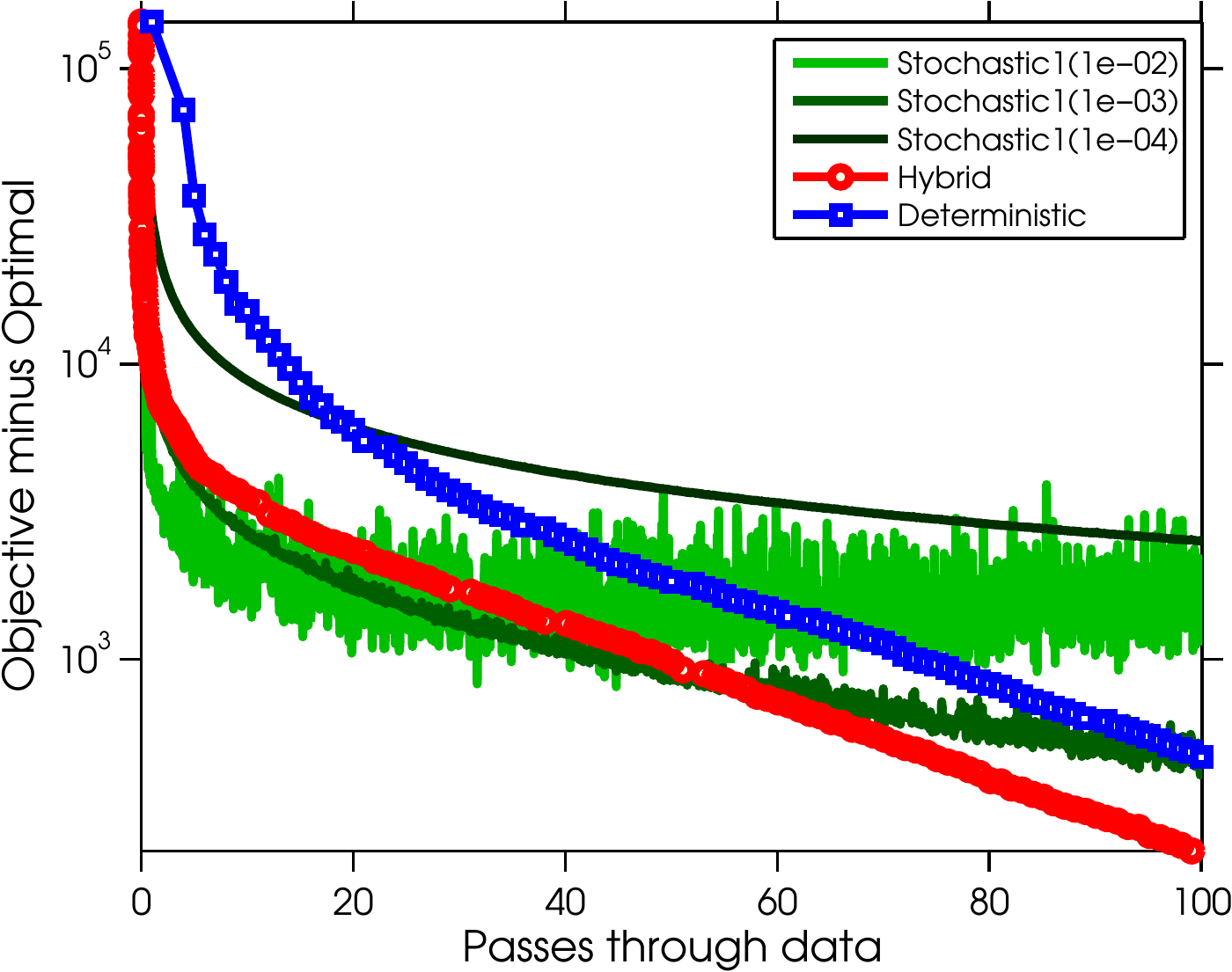}
  \hfill
  \includegraphics[width=.48\textwidth]{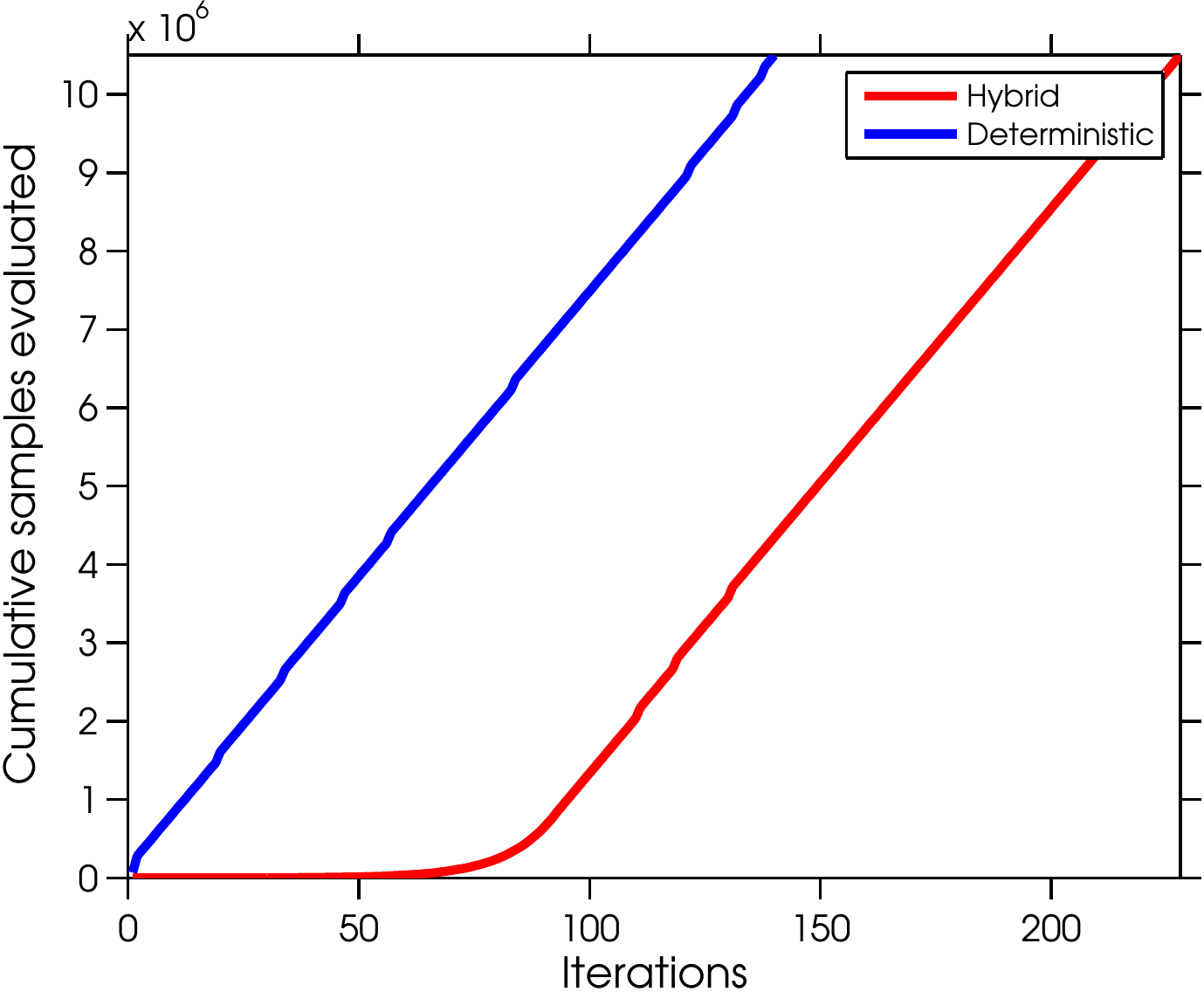}
  \caption{Multinomial logistic regression experiments for different optimization
    strategies for digit classification.}
  \label{fig:mlogreg}
\end{figure}
 
Our experiments for multinomial logistic regression are based on the
well-known MNIST data set~\cite{mnist}, containing 70,000 examples of
28-by-28 images of digits, where each digit is classified as one of
the numbers 0 through 9. The results are plotted in
Figure~\ref{fig:mlogreg}, with $\lambda=1$. The trends are similar to
the binary logistic regression experiments.

\subsection{Chain-structured conditional random fields}
\label{sec:crf}

The binary and multinomial logistic regression models consider the
case where a particular outcome $b_i$ is associated with a particular
feature vector $a_i$.  The CRF model takes multinomial logistic
regression a step further by considering a set of feature vectors
$a_i^k$ with which to predict corresponding values for discrete
variables $b_i^k \in \Cscr$, where $k$ takes values in a discrete
ordered set $\Omega$.

We can naturally extend the multinomial logistic model to this
scenario by defining the probability of a joint assignment $b_i$
as
\[
p(\{b_i^k = j_k\}_{k\in\Omega} \ \mid \{a_i^k\}_{k\in\Omega},
  \{x_j\}_{j\in\Cscr}) =
\frac{1}{Z_i}\prod_{k\in\Omega} \exp(x_{j_k}^T a_i^k).
\]
(We assume that the parameter vectors $x_{j_k}$ are tied, so that
$x_{j_k}$ is constant for all $k$; this assumption is not required in
general.)  The normalizing constant
\[
Z_i = \sum_{j_1\in\Cscr}\sum_{j_2\in\Cscr}\dots\sum_{j_k\in\Cscr}
      \prod_{k\in\Omega}\exp(x_{j_k}^T a_i^k)
\]
is chosen so that the distribution sums to one over all possible
configurations of the $b_i^k$ variables.  As written, computing $Z_i$ involves a very
large number of terms, $\abs{\Cscr}^{\abs{\Omega}}$; still, the sum can
be computed efficiently by exchanging the order of operations.  While
this model is a straight-forward generalization of multinomial
logistic regression, it assumes that the labels in $\Omega$ that we
are simultaneously predicting are independent.  This might be
unrealistic if, for example, the variables come from time-series data
where $b_i^k$ and $b_i^{k+1}$ are likely to be correlated.

\begin{figure}[t]
  \centering
  \includegraphics[width=.48\textwidth]{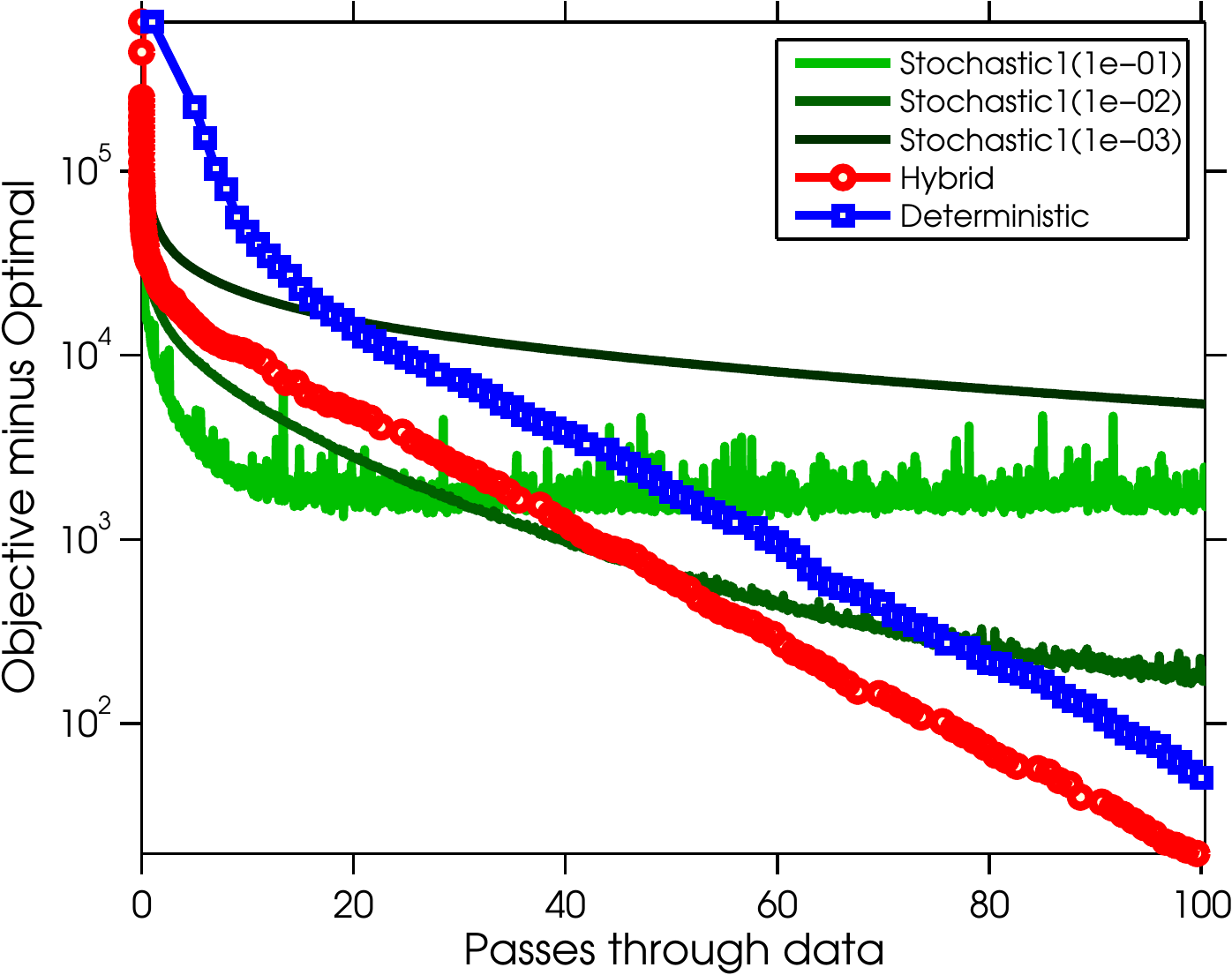}
  \hfill
  \includegraphics[width=.48\textwidth]{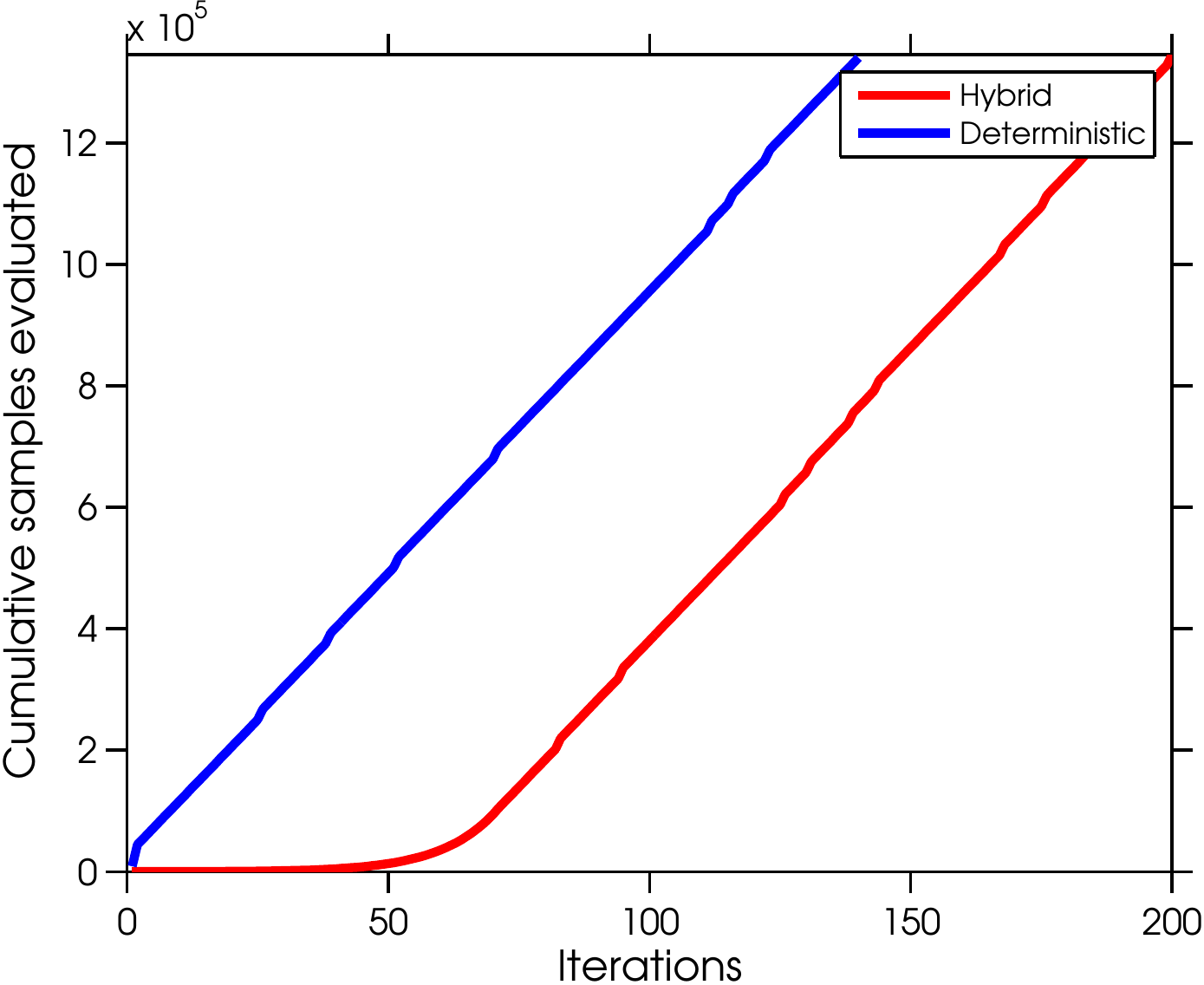}
  \caption{Chain-structured conditional random field experiments for
    different optimization strategies on the noun-phrase chunking
    task.}
  \label{fig:crf}
\end{figure}
 
A chain-structured CRF~\cite{lafferty2001conditional} augments the
model with additional terms that take into account sequential
dependencies in the labels.  It allows pairwise features $a_i^{kk'}$
and associated parameters $x_{kk'}$, and defines the probability of a
configuration as
\begin{multline*}
  p_3(\{b^k_i = j\}_{k\in\Omega} \mid \{a_i^k\}_{k\in\Omega}, \{a_i^{kk'}\}_{k,k'\in\Omega,k'=k+1},
  \{x_j\}_{j\in\Cscr}, \{x_{j,j'}\}_{j,j'\in\Cscr})
  \\= \frac{1}{Z_i}
    \left[\prod_{\vphantom{k'\atop k'}k\in\Omega}
      \exp(x_{j_k}^T a_i^k)
    \right]
    \cdot
    \left[
      \prod_{k,k'\in\Omega\atop k'=k+1}\exp(x_{j_kj_{k'}}^T a_i^{kk'})
    \right]
    .
\end{multline*}
The normalizing constant $Z_i$ is again set so that the distribution
sums to one, and it can be computed using a variant on the
forward-backward algorithm used in hidden Markov
models~\cite[\S~III]{rabiner1989tutorial}.
Because we need to run the forward-backward
algorithm for each $i$, the probability function $p_3$ is
significantly more expensive to evaluate than the corresponding
multinomial probability function $p_2$;
see~\eqref{eq:19}.

For our chain-structured CRF experiments, we use the noun-phrase
chunking problem from the CoNLL-2000 Shared Task~\cite{SangBuc:2000},
where the goal is to assign each word in a sentence to one of 22
possible states, and we use approximately 1.6 million features to
represent the presence of words and word
combinations~\cite{SangBuc:2000,sha2003shallow}.  The CoNLL-2000 data
set contains 211,727 words grouped into 8,936 sentences.  The results
of this experiment are plotted in Figure~\ref{fig:crf}, where the
regularization parameter is $\lambda=1$.

\subsection{General conditional random fields}
\label{sec:g-crf}

\begin{figure}[t]
  \centering
  \includegraphics[width=.48\textwidth]{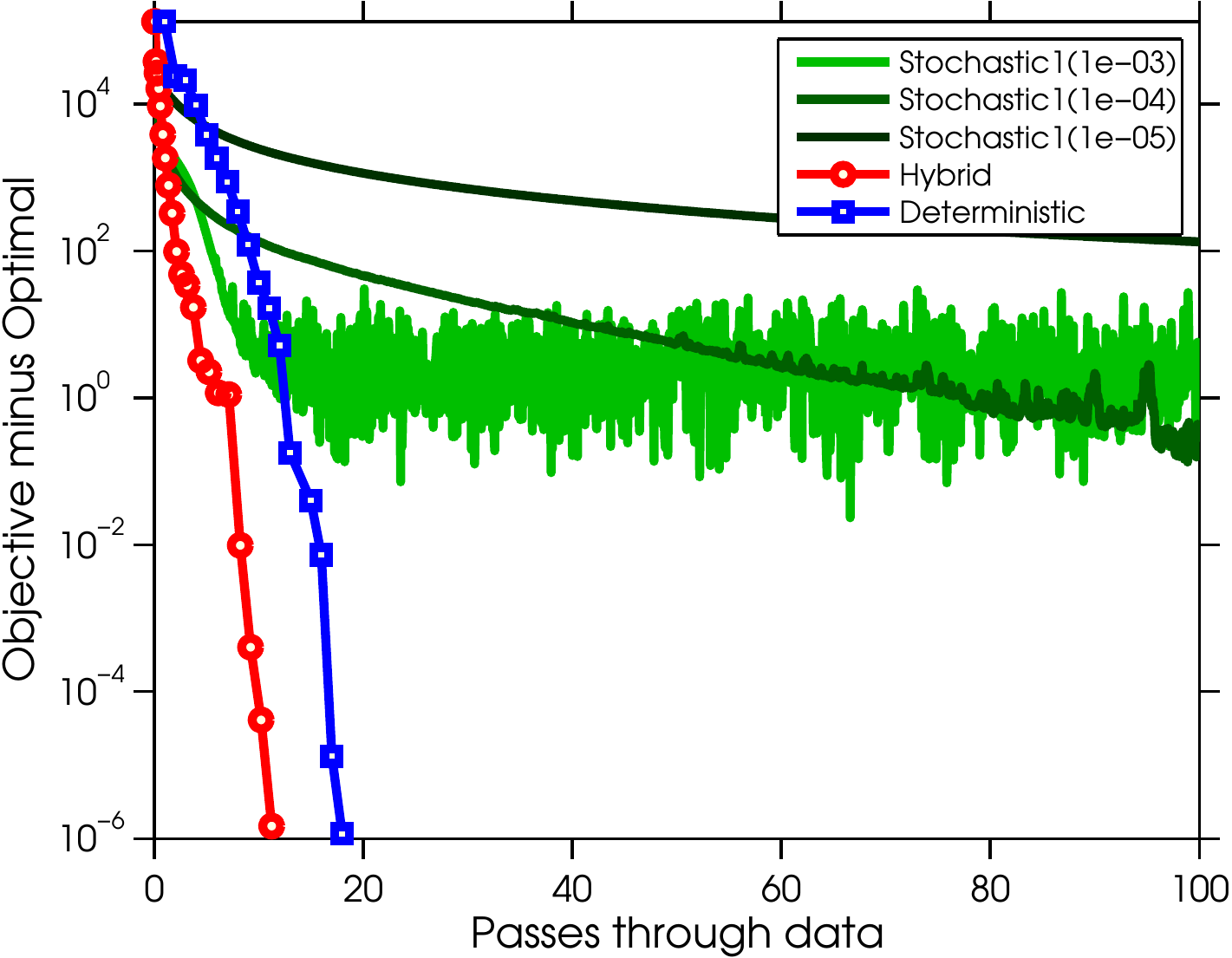}
  \hfill
  \includegraphics[width=.48\textwidth]{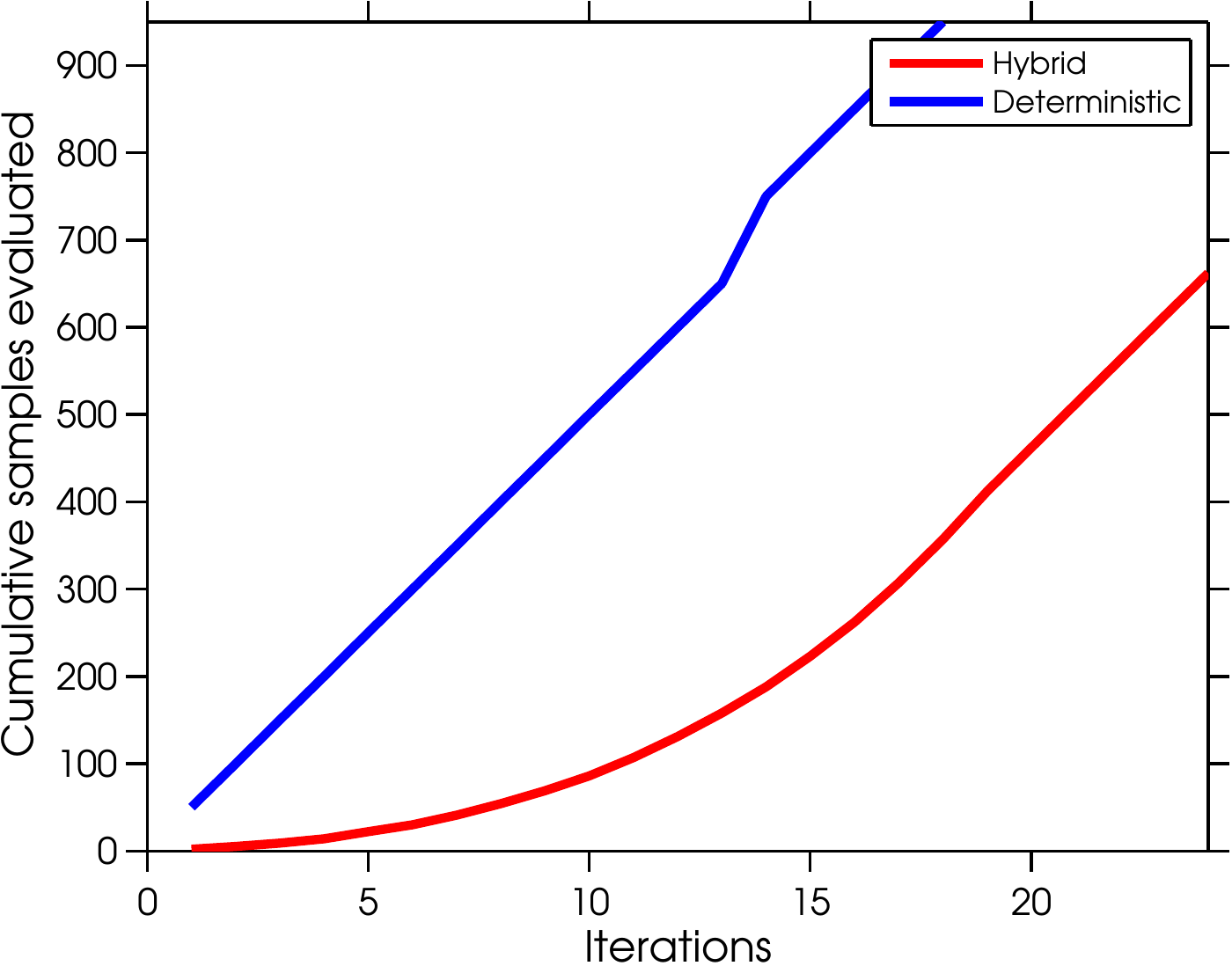}
  \caption{Lattice-structured conditional random field experiments for
    different optimization strategies for an image denoising task.}
  \label{fig:crf2d}
\end{figure}
 
While chain-structures can model sequential dependencies in the
labels, 
we might be interested in other structures, such
as lattice structures for modeling image data.  In order to use
general CRFs~\cite[\S4.6,1 and \S20.3.2]{koller2009pgm}, we
define a graph $\Gscr$ where the nodes are the labels
$\{1,2,\dots,w\}$ and the edges are the variable dependencies that
we wish to consider.  We then define the probability of $b_i$ taking a
configuration $b_i^k=j_k$, for $k\in\Omega$, as
\begin{multline*}
  p_4(\{b_i^k=j_k\}_{k\in\Omega} \mid \{a_i^{k'}\}_{k'\in\Omega},\{a_i^{kk'}\}_{(k,k')\in\Escr},
    \{x_j\}_{j\in\Cscr},\{x_{jj'}\}_{(j,j')\in\Cscr})
    \\= \frac{1}{Z_i}
    \left[
      \prod_{k\vphantom'\in\Omega} \exp(x_{j_k}^T a_i^k)
    \right]
    \cdot
    \left[
      \prod_{(k,k')\in\Escr}\exp(x_{j_kj_{k'}}^T a_i^{kk'})
    \right]
    .
\end{multline*}

\begin{figure}[t] 
  \centering\small\setlength{\fboxsep}{1pt}
  \begin{tabular}{@{}c@{\hspace{.1in}}c@{}}
        \fbox{\includegraphics[width=.30\textwidth]{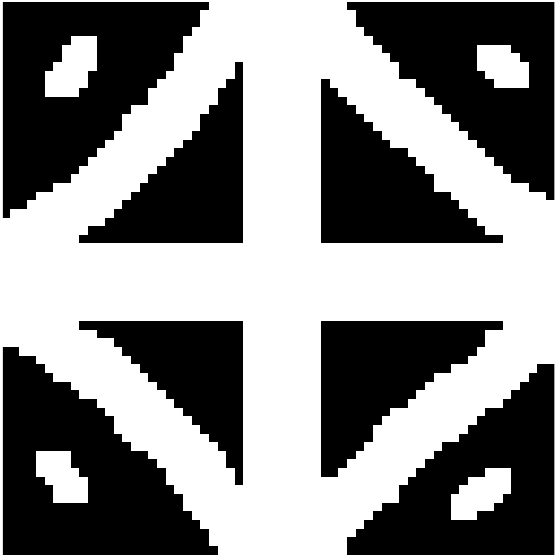}}
       &\fbox{\includegraphics[width=.30\textwidth]{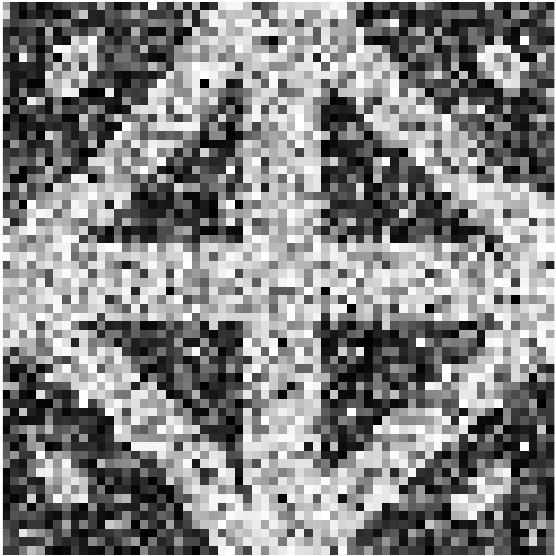}}
     \\(a) original image & (b) noisy image
  \end{tabular}
  \\[3pt]
  \begin{tabular}{@{}c@{\hspace{.1in}}c@{\hspace{.1in}}c@{}}
     \fbox{\includegraphics[width=.30\textwidth]{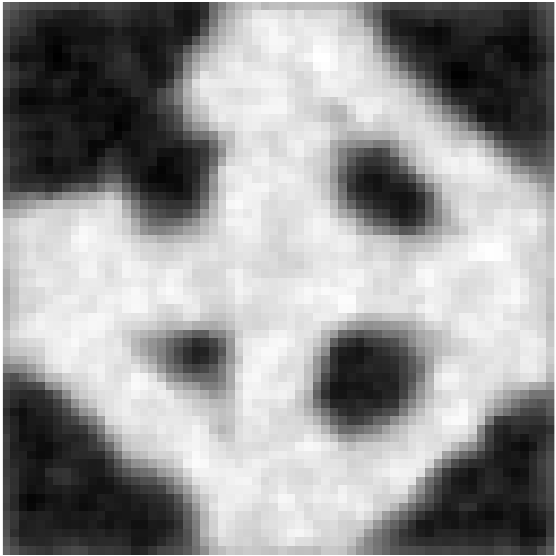}}
   & \fbox{\includegraphics[width=.30\textwidth]{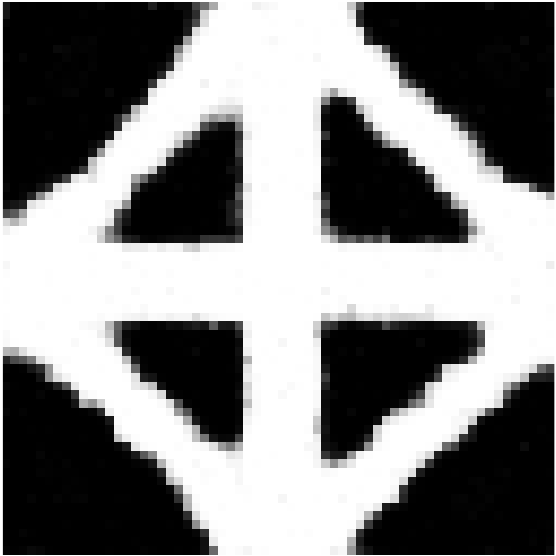}}
   & \fbox{\includegraphics[width=.30\textwidth]{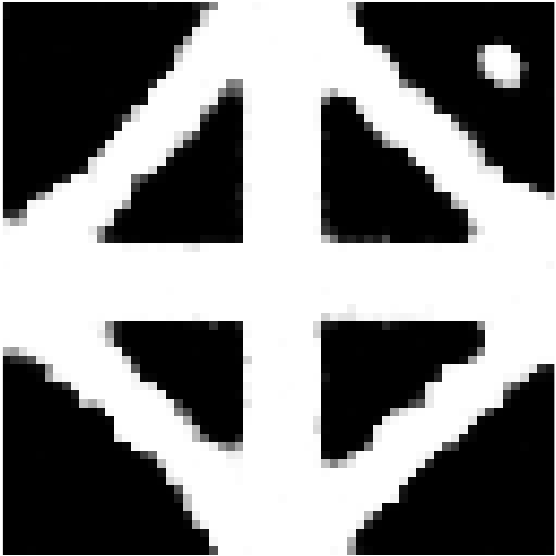}}
  \\ (c) & (d) & (e)
  \\[3pt]
     \fbox{\includegraphics[width=.30\textwidth]{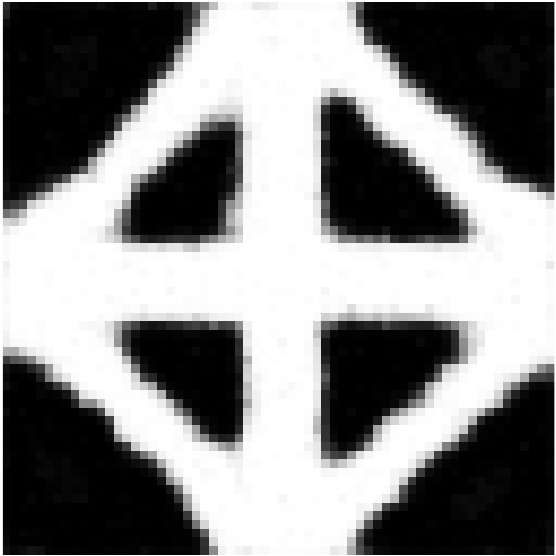}}
   & \fbox{\includegraphics[width=.30\textwidth]{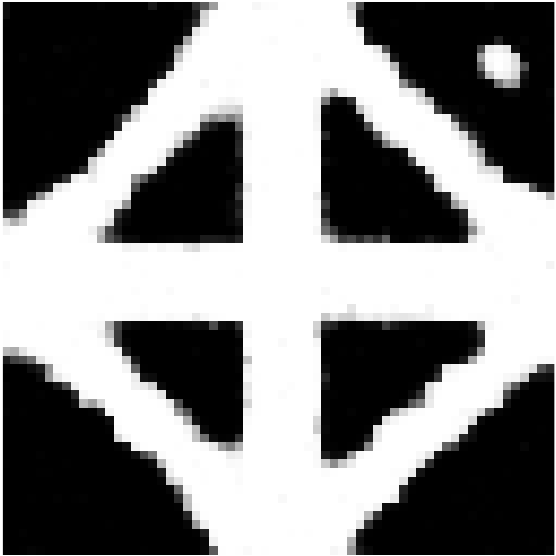}}
   & \fbox{\includegraphics[width=.30\textwidth]{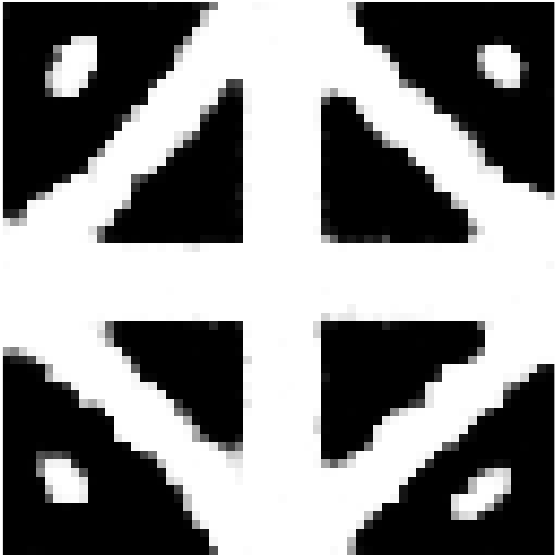}}
  \\ (f) & (g) & (h)
  \end{tabular}
  \caption{Top row: original (a) and noisy (b) image.  
    Second
    row: marginals after 2 passes through the data for deterministic
    (c), stochastic (d), and hybrid (e).  Third row:
    marginals after 5 passes through the data for deterministic
    (f), stochastic (g), and hybrid (h).}
  \label{fig:crf2d_images}
\end{figure}
 
In this general case, computing $Z_i$ is in the complexity
class $\sharp P$, and the best known algorithms have a runtime that is
exponential in the tree-width of the graph
$\Gscr$~\cite[\S9-10]{koller2009pgm}.  For a two-dimensional
lattice structure, the tree-width is the minimum between the two
dimensions of the structure, so computing $Z_i$ is only feasible if
one of the dimensions is very small (in the
degenerate one-dimensional chain-structured case, the tree-width is
one).

Because of the intractability of computing $Z_i$, we consider
optimizing a pseudo-likelihood
approximation~\cite{besag1975statistical} based on the probability
model
\begin{multline*}
  p_4(\{b_i^k=j_k\}_{k\in\Omega} \mid \{a_i^k\}_{k\in\Omega},\{a_i^{kk'}\}_{(k,k')\in\Escr},
      \{x_j\}_{j\in\Cscr},\{x_{jj'}\}_{(j,j')\in\Cscr})
 \\ \approx \prod_{k\in\Omega}
     p(b_i^k = j_k \mid \{a_i^{k'}\}_{k'\in\Omega}, \{a_i^{kk'}\}_{(k,k')\in\Escr}, \{x_j\}_{j\in\Cscr},\{x_{jj'}\}_{(j,j')\in\Cscr},
     \{b_i^{k'}\}_{k'\in\Omega,k'\ne k}).
\end{multline*}
The individual terms in this product of conditionals have the form of
a multinomial logistic regression probability and are straightforward
to compute. 
This is the function $p$
used to define the objective function $f$ in~\eqref{eq:18}.

Our experiments on general CRFs are based on the image-denoising
experiments described by Kumar and
Hebert~\cite{kumar2004discriminative}. We use their set of 50
synthetic 64-by-64 images.  Figure~\ref{fig:crf2d} shows the
performance of the different methods with a regularization parameter of
$\lambda=1$.  Figure~\ref{fig:crf2d_images} illustrates the marginal
probabilities for the different methods at various points in the
optimization for a randomly-chosen image in the data set. (For the
stochastic method, we plot the result with a step size of
$\alpha=10^{-4}$.)  To approximate these marginals, we use the loopy
belief propagation message-passing
algorithm~\cite[\S8.4.7]{bishop2006pattern}.  In these plots we see
that the deterministic method does poorly even after two full passes
through the data set, while the stochastic and hybrid methods do much
better.  After five passes through the data set, the hybrid method has
found a solution that is visually nearly indistinguishable from the
true solution, while it is still possible to see obvious differences
in the deterministic and stochastic methods.

\subsection{Seismic inversion}
\label{sec:nnls}

This last numerical experiment is a seismic inversion problem
described by van Leeuwen et al.~\cite{Tristan11}. The aim here is to
recover an image of underground geological structures using only data
collected by geophone receivers placed at the surface of the earth; these
geophones record acoustic ``shots'' created by sources also at the
surface.

\begin{figure}[t]
  \centering
  \includegraphics[width=.48\textwidth]{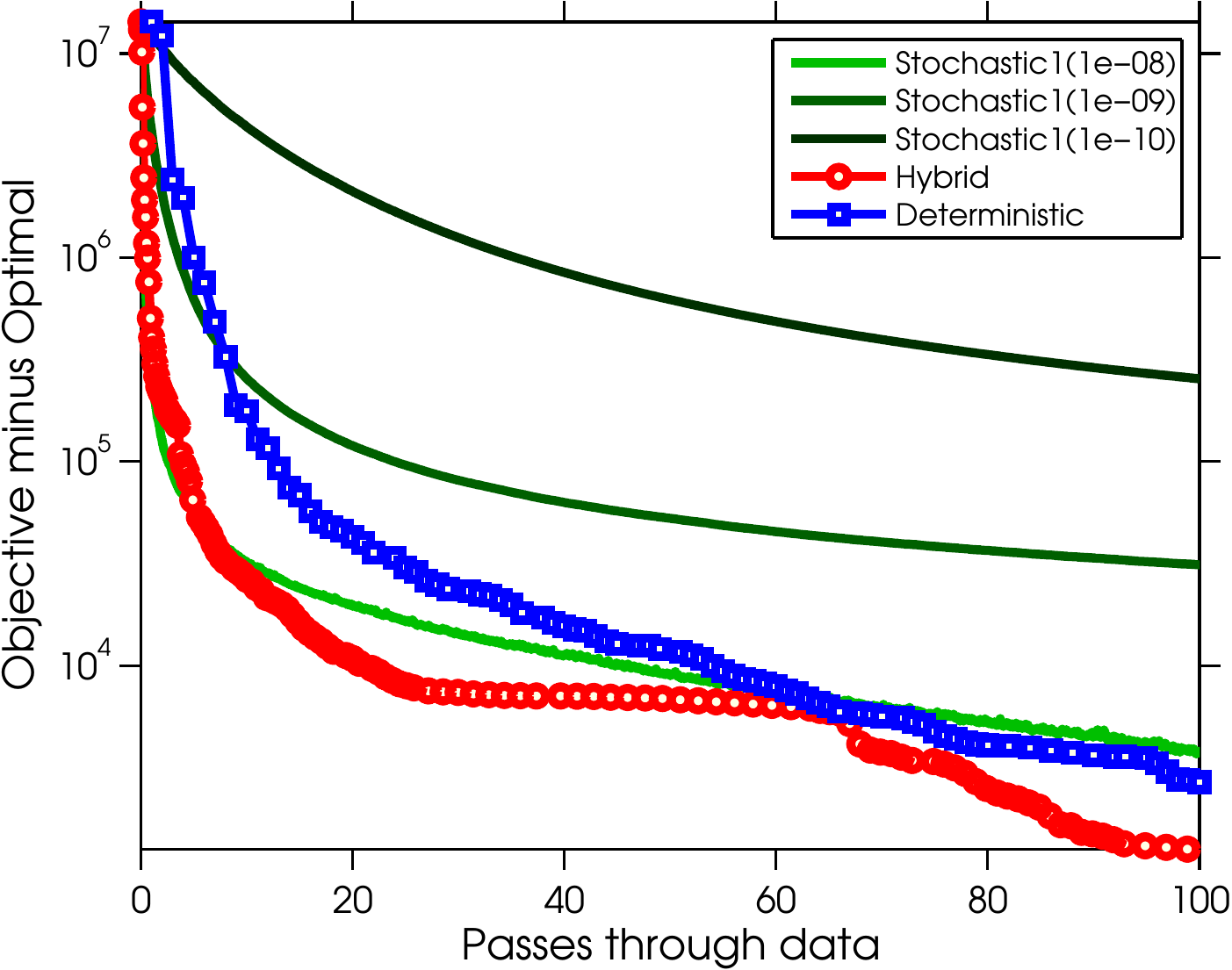}
  \hfill
  \includegraphics[width=.48\textwidth]{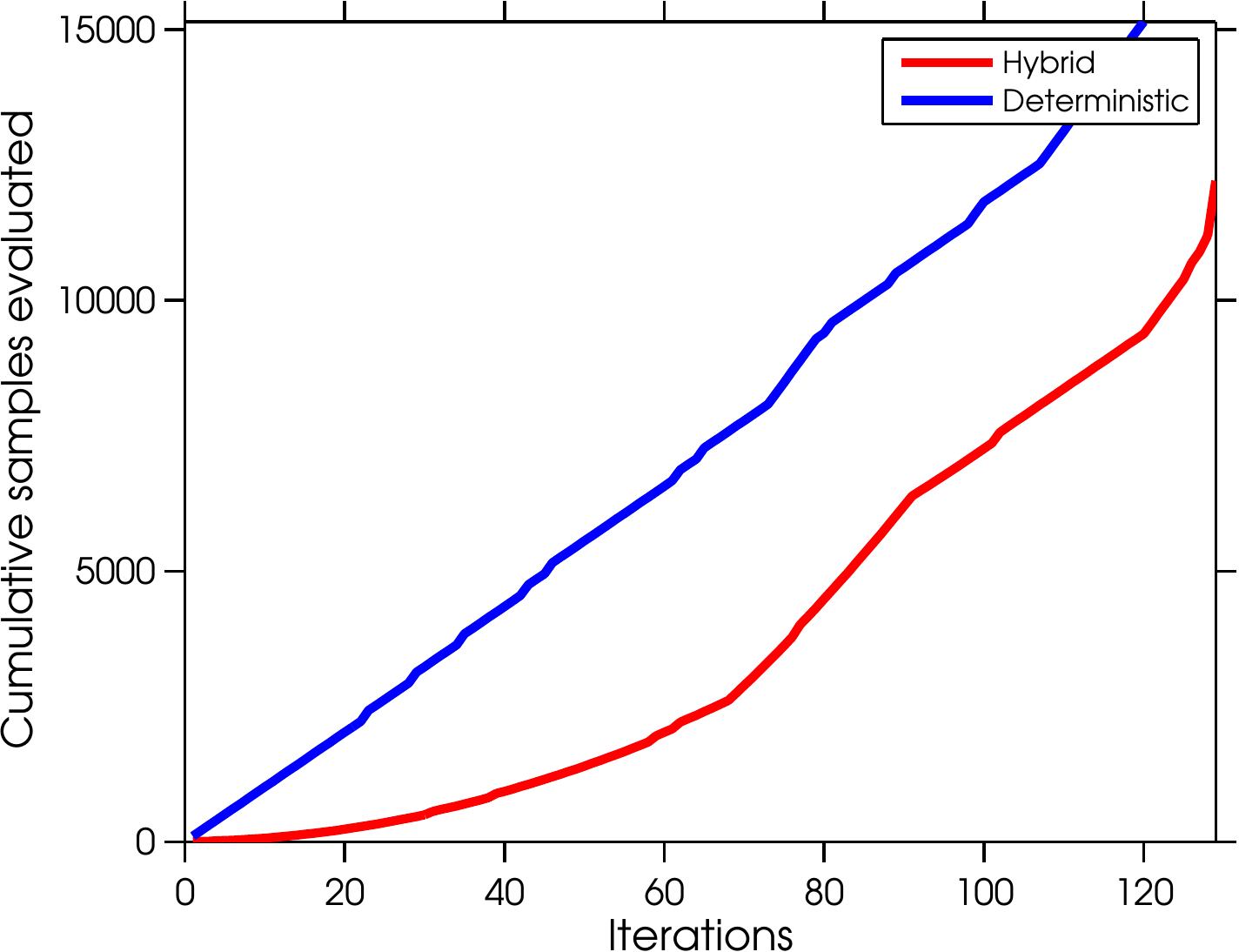}
  \caption{Nonlinear least-squares experiments for different optimization
    strategies on a seismic inversion problem.}
  \label{fig:dwi}
\end{figure}
 
The waveform inversion problem attempts to find a model $x$ of the
subsurface structure that minimizes the nonlinear least-squares misfit
as measured by the function
\[
 \phi(x) =
 \sum_{i=1}^M \sum_{\omega\in\Omega}\norm{d_i - PH_\omega(x)\inv q_i}^2.
\]
Each index $i$ corresponds to a particular shot (i.e., an observation)
created by the source $q_i$ which creates a measurement $d_i$; each
experiment samples a set of frequencies $\omega\in\Omega$. The matrix
$P$ samples the wavefield at the receiver locations. The main cost in
evaluating $\phi$ is solving the Helmholtz equation $H_\omega[x] u =
q_i$ for each $i$, which is an expensive partial differential
equation. Regularization is typically achieved by truncating the
solution process~\cite{Tristan11}. Although the function $\phi$ is
nonconvex and does not satisfy our standing assumptions, it hints at
the applicability of the hybrid approach for solving difficult
problems of important practical interest.

The results shown in Figure~\ref{fig:dwi} are based on a relatively
small 2-dimensional example that involves M=101 sources measuring 8
frequencies. (Larger experiments, especially in 3-dimensions, are only
feasibly accomplished on a computing cluster.)  
Because the problem is not strongly convex and hence we do not
expect a fast convergence rate when far from this solution, for this problem the hybrid
method increases the sample size by only one sample at each iteration, i.e.,
$\abs{\batch\kp1}=\min\{\abs{\batch\k}+1,\, M\}$. After 20 passes
through the data, the hybrid method clearly outperforms the
deterministic method. The best stochastic method performs nearly as
well as the deterministic and hybrid methods for the first 60
iterations, though with other step sizes it performs poorly. Although
the figure shows the best methods all achieving similar residuals
after 60 passes through the data, in practice that involves a
prohibitive number of Helmholtz solves; practitioners are interested
in making quick progress with as few solves as possible (and without
needing to test a variety of step sizes to achieve good performance).

\section{Discussion}\label{sec:discussion}

Our work has focused on inexact gradient methods for the unconstrained
optimization of differentiable strongly-convex objectives. We
anticipate that a similar convergence analysis with inexact gradients
could be applied to other algorithms, such as Nesterov's accelerated
gradient method~\cite[\S2.1]{nesterov2004introductory}. In this case,
it may be possible to relax the strong convexity assumption, and
obtain the optimal $\Oscr(1/k^2)$ rate for an inexact-gradient version
of this algorithm.  The optimal $\Oscr(1/k^2)$ rate using a gradient
approximation whose error is uniformly bounded across iterations has
been established by several authors, notably
d'Aspremont~\cite{d2008smooth}. But allowing a variable error would
encourage more flexibility in the early iterations, and would allow
for eventually solving the problem to arbitrary accuracy.  Although
the $\Oscr(1/k^2)$ rate is sub-linear, it is substantially faster than
the optimal $\Oscr(1/\sqrt{k})$ achievable by methods that use noisy
gradient information~\cite[\S14.1]{nemirovski1994efficient}.

We might also obtain analogous rates for proximal-gradient methods for
optimization with convex constraints or non-differentiable composite
optimization problems, such as 1-norm
regularization~\cite{nesterov2007gradient}.  The more general class of
mirror-descent methods~\cite{beck2003mirror}, which are useful for
problems with a certain geometry such as optimization with simplex
constraints, also seem amenable to analysis in our controlled
error scenario.

We considered the case of bounded noise or noise that can be bounded
in expectation, and subsequently derived convergence rates and
expected convergence rates, respectively. We might also consider the
case where the noise is bounded with a certain probability.  If the
individual $\nabla f_i(x)$ are concentrated around $\nabla f(x)$, this
might allow us to use concentration
inequalities~\cite{massart2007concentration} to show that the
convergence rates hold with high probability.  Although we have
analyzed an arbitrary strategy for selecting the elements of the
sample, and shown that uniform sampling achieves a better bound in
expectation, it is possible that a quasi-random selection of the
individual gradients might further refine the
bound~\cite{morokoff1994quasi}.


Although our emphasis here is on data fitting applications where the
error is a by-product of subsampling the data, our analysis and
implementation may be useful for other problems. For example, Gill et
al.~\cite[p.~357]{gill1981practical} discuss the case of an objective
function that can be evaluated to a prescribed accuracy (e.g., it
could depend on an iterative process or a discretization level).  They
suggest solving the optimization problem over a sequence of tighter
function accuracies. Our work provides a formal analysis and practical
implementation of a method where the accuracy might be increased
dynamically each iteration rather than solving a sequence of
intermediate optimization problems.  As a more recent example,
Poyiadjis et al.~\cite{poyiadjis2011particle} consider approximating
gradients in non-Gaussian state-space models using particle filters.
Here, the variance of the approximation is directly proportional to
the number of particles, and thus our work provides guidelines for
selecting the number of particles to use in the approximation on each
iteration.

\section*{Acknowledgements}

Thanks to Kevin Swersky for the suggesting the growing sample-size
strategy, and to Sasha Aravkin, Francis Bach, Hoyt Koepke, Pierre
Jacob, Simon Lacoste-Julien, and Nicolas Le Roux for valuable
discussions. We are grateful to two anonymous referees for their
thoughtful comments, which led to many detailed changes, including the
addition of \S\ref{sec:weakly}. We thank Ming Yan for highlighting a
mistake in the original proof of Theorem~\ref{th:incremental-weak}.

\bibliographystyle{siam}
\bibliography{bib}

\end{document}